\documentclass[nonblindrev]{informs2017}
\OneAndAHalfSpacedXI 

\usepackage{endnotes}
\let\footnote=\endnote
\let\enotesize=\normalsize
\def\notesname{Endnotes}%
\def\makeenmark{\hbox to1.275em{\theenmark.\enskip\hss}}
\def\enoteformat{\rightskip0pt\leftskip0pt\parindent=1.275em
\leavevmode\llap{\makeenmark}}

\usepackage{multirow}
\usepackage{float}
\usepackage{epstopdf}
\usepackage{appendix}
\usepackage{bm}
\usepackage{dsfont}

\usepackage{graphicx}
\usepackage{subcaption}
\usepackage{caption}
\usepackage{enumitem}
\usepackage{array}

\newcommand{\tabincell}[2]{\begin{tabular}{@{}#1@{}}#2\end{tabular}}

\newcommand{\bmt}[1]{\tilde{\bm{#1}}}
\newcommand{\bmh}[1]{\hat{\bm{#1}}}

\usepackage{tabularx}
\usepackage{pict2e}
\newcolumntype{L}[1]{>{\raggedright\arraybackslash}p{#1}}
\newcolumntype{C}[1]{>{\centering\arraybackslash}p{#1}}
\newcolumntype{R}[1]{>{\raggedleft\arraybackslash}p{#1}}

\usepackage{color}
\usepackage[usenames,dvipsnames]{xcolor}
\definecolor{strcolor}{rgb}{0.6, 0.2, 0.6}
\definecolor{commentcolor}{rgb}{0.3125, 0.5, 0.3125}
\definecolor{keycol}{rgb}{0, 0, 1}

\usepackage{listings}
\lstdefinelanguage{Julia}%
{morekeywords={abstract,break,case,catch,const,continue,do,else,elseif,%
		end,export,false,for,function,immutable,import,importall,if,in,%
		macro,module,otherwise,quote,return,switch,true,try,type,typealias,%
		using,while},%
sensitive=true,%
alsoother={\$},%
morecomment=[l]\#,%
morecomment=[n]{\#=}{=\#},%
morestring=[s]{"}{"},%
morestring=[m]{'}{'},%
}[keywords,comments,strings]%

\usepackage{natbib}
 \bibpunct[, ]{(}{)}{,}{a}{}{,}%
 \def\bibfont{\small}%
\TheoremsNumberedThrough     
\ECRepeatTheorems

\EquationsNumberedThrough    

\begin{document}


\RUNAUTHOR{Chen, Kuhn, and Wiesemann}

\RUNTITLE{Data-Driven Chance Constrained Programs over Wasserstein Balls}

\TITLE{Data-Driven Chance Constrained Programs \\ over Wasserstein Balls}

\ARTICLEAUTHORS{%
\AUTHOR{Zhi Chen}
\AFF{College of Business, City University of Hong Kong, Kowloon Tong, Hong Kong, \\ zhi.chen@cityu.edu.hk}
\AUTHOR{Daniel Kuhn}
\AFF{Risk Analytics and Optimization Chair, \'{E}cole Polytechnique F\'{e}d\'{e}rale de Lausanne, Lausanne, Switzerland, \\ daniel.kuhn@epfl.ch}
\AUTHOR{Wolfram Wiesemann}
\AFF{Imperial College Business School, Imperial College London, London, United Kingdom, \\ ww@imperial.ac.uk}
}
\ABSTRACT{We provide an exact deterministic reformulation for data-driven chance constrained programs over Wasserstein balls. For individual chance constraints as well as joint chance constraints with right-hand side uncertainty, our reformulation amounts to a mixed-integer conic program. In the special case of a Wasserstein ball with the $1$-norm or the $\infty$-norm, the cone is the nonnegative orthant, and the chance constrained program can be reformulated as a mixed-integer linear program. Our reformulation compares favourably to several state-of-the-art data-driven optimization schemes in our numerical experiments. 
}%


\KEYWORDS{Distributionally robust optimization; ambiguous chance constraints; Wasserstein distance.}

\HISTORY{\today}

\maketitle

%


\section{Introduction}

Distributionally robust optimization is a powerful modeling paradigm for optimization under uncertainty, where the distribution of the uncertain problem parameters is itself uncertain, and where the performance of a decision is assessed in view of the worst-case distribution from a prescribed ambiguity set. The earlier literature on distributionally robust optimization has focused on moment ambiguity sets which contain all distributions that obey certain (standard or generalized) moment conditions; see, {\em e.g.}, \citet{Delage_Ye_2010}, \citet{Goh_Sim_2010} and \citet{Wiesemann_Kuhn_Sim_2014}. \citet{Pflug_Wozabal_2007} were the first to propose an ambiguity set of the form of a ball in the space of distributions with respect to the celebrated Wasserstein, Kanthorovich or optimal transport distance. The type-1 Wasserstein distance $d_{\rm W}(\mathbb{P}_1,\mathbb{P}_2)$ between two distributions $\mathbb{P}_1$ and $\mathbb{P}_2$ on $\mathbb{R}^K$, equipped with a general norm $\|\cdot\|$, is defined as the minimal transportation cost of moving $\mathbb{P}_1$ to $\mathbb{P}_2$ under the premise that the cost of moving a Dirac point mass from $\bm\xi_1$ to $\bm\xi_2$ amounts to $\|\bm{\xi}_1 - \bm{\xi}_2\|$. Mathematically, this implies that
$$
\begin{array}{rcl}
d_{\rm W}(\mathbb{P}_1,\mathbb{P}_2) \; = \; &\displaystyle \inf_{\mathbb{P} \in \mathcal{P}(\mathbb{P}_1, \mathbb{P}_2)}
& \mathbb{E}_{\mathbb{P}}[\|\bmt{\xi}_1 - \bmt{\xi}_2\|] ,
\end{array}
$$
where $\bmt{\xi}_1 \sim \mathbb{P}_1, \bmt{\xi}_2 \sim \mathbb{P}_2$, and $\mathcal{P}(\mathbb{P}_1, \mathbb{P}_2)$ represents the set of all distributions on $\mathbb{R}^K\times \mathbb{R}^K$ with marginals $\mathbb{P}_1$ and $\mathbb{P}_2$. The Wasserstein ambiguity set $\mathcal{F}(\theta)$ is then defined as a ball of radius $\theta\ge 0$ with respect to the Wasserstein distance, centered at a prescribed reference distribution~$\hat{\mathbb{P}}$:
\begin{equation}\label{set:Wasserstein}
\mathcal{F}(\theta) = \{\mathbb{P} \in \mathcal{P}(\mathbb{R}^K) \mid
d_{\rm W}(\mathbb{P}, \hat{\mathbb{P}}) \leq \theta\}.
\end{equation}
One can think of the Wasserstein radius $\theta$ as a budget on the transportation cost. Indeed, any member distribution in $\mathcal{F}(\theta)$ can be obtained by rearranging the reference  distribution $\hat{\mathbb{P}}$ at a transportation cost of at most $\theta$. If only a finite training dataset $\{\bmh{\xi}_i\}_{i \in [N]}$ is available, a natural choice for $\hat{\mathbb{P}}$ is the empirical distribution $\hat{\mathbb{P}} = \frac{1}{N}\sum_{i = 1}^N \delta_{\bmh{\xi}_i}$, which represents the uniform distribution on the training samples. Throughout the paper, we will assume that $\hat{\mathbb{P}}$ is the empirical distribution.

While it has been recognized early on that Wasserstein ambiguity sets offer many conceptual advantages (\emph{e.g.}, their member distributions do not need to be absolutely continuous with respect to $\hat{\mathbb{P}}$ and, if properly calibrated, they constitute confidence regions for the unknown true data-generating distribution), it was believed that they almost invariably lead to hard global optimization problems. Recently, \citet{Esfahani_Kuhn_2017} and \citet{Zhao_Guan_2018} discovered that many interesting distributionally robust optimization problems over Wasserstein ambiguity sets can actually be reformulated as tractable convex programs---provided that $\hat{\mathbb{P}}$ is discrete and that the problem's objective function satisfies certain convexity properties. These reformulations have subsequently been generalized to Polish spaces and non-discrete reference distributions by \citet{blanchet2019quantifying} and \citet{Gao_Kleywegt_2016}. Since then, distributionally robust optimization models over Wasserstein ambiguity sets have been proposed for many applications, including transportation (\citealt{carlsson2018wasserstein}) and machine learning (\citealt{blanchet2019robust}, \citealt{gao2017distributional}, \citealt{shafieezadeh2019regularization} and \citealt{sinha2017certifiable}).

In this paper we study distributionally robust chance constrained programs of the form
\begin{equation}\label{prob:cc general}
\begin{array}{cll}
\displaystyle \min_{\bm{x} \in \mathcal{X}} &~\bm{c}^\top\bm{x} \\
{\rm s.t.} &~\displaystyle \mathbb{P}[\bmt{\xi} \in \mathcal{S}(\bm{x})] \geq 1-\varepsilon &~\forall \mathbb{P} \in \mathcal{F}(\theta),
\end{array}
\end{equation}
where the goal is to find a decision $\bm{x}$ from within a compact polyhedron $\mathcal{X} \subseteq \mathbb{R}^L$ that minimizes a linear cost function $\bm{c}^\top\bm{x}$ and ensures that the exogenous random vector $\bmt{\xi}$ falls within a decision-dependent safety set $\mathcal{S}(\bm{x}) \subseteq \mathbb{R}^K$ with high probability $1-\varepsilon$ under every distribution $\mathbb{P} \in \mathcal{F}(\theta)$. Since the reference distribution $\hat{\mathbb{P}}$ in~\eqref{prob:cc general} is the empirical distribution over the training dataset $\{\bmh{\xi}_i\}_{i \in [N]}$, we refer to~\eqref{prob:cc general} as a \emph{data-driven} chance constrained program.

To date, the literature on data-driven chance constraints has focused primarily on variants of problem~\eqref{prob:cc general} where the Wasserstein ambiguity set $\mathcal{F} (\theta)$ is replaced with an ambiguity set $\mathcal{G} (\theta)$ that contains all distributions close to the empirical distribution $\hat{\mathbb{P}}$ with respect to a $\phi$-divergence (such as the Kullback-Leibler divergence or the $\chi^2$-distance):
\begin{equation*}
\mathcal{G}(\theta) = \bigg\{\mathbb{P} \in \mathcal{P}(\mathbb{R}^K) ~\bigg|~\mathbb{P}\ll \hat{\mathbb{P}}, \;\;
\int_{\mathbb{R}^K} \phi\bigg(\dfrac{{\rm d}\mathbb{P}(\bm{\xi})}{{\rm d}\hat{\mathbb{P}}(\bm{\xi})}\bigg){\rm d}\hat{\mathbb{P}}(\bm{\xi}) \leq \theta \bigg.\bigg\},
\end{equation*}
where $\phi: \mathbb{R}_+ \rightarrow \mathbb{R}$ is the divergence function. \cite{Hu_Hong_2013} show that a distributionally robust chance constrained program over a Kullback-Leibler ambiguity set reduces to a classical chance constrained progam over the reference distribution $\hat{\mathbb{P}}$ and an adjusted risk threshold $\varepsilon' < \varepsilon$. While this result holds for any reference distribution, $\phi$-divergence ambiguity sets only contain distributions that are absolutely continuous with respect to $\hat{\mathbb{P}}$, that is, any distribution in $\mathcal{G} (\theta)$ only assigns positive probability to those measurable subsets $A \subseteq \mathbb{R}^K$ for which $\hat{\mathbb{P}} [\tilde{\bm{\xi}} \in A] > 0$. This is undesirable for problems with a large dimension $K$ and/or few training data, where it is unlikely that every possible value of $\tilde{\bm{\xi}}$ has been observed in $\{\bmh{\xi}_i\}_{i \in [N]}$. This shortcoming is addressed by \cite{Jiang_Guan_2016, jiang2018risk}, who replace the reference distribution with a Kernel density estimator.

Despite their tremendous success and widespread adoption in recent years, the use of $\phi$-divergences can lead to undesirable side effects in some applications: they compare distributions on a ``scenario-by-scenario" basis and thus do not consider the possibility of noisy measurements \citep{Gao_Kleywegt_2016}, and they generically fail to be probability metrics as they typically violate symmetry as well as the triangle inequality. Moreover, as we show next, $\phi$-divergence ambiguity sets may be overly optimistic when only few training samples are available.

\vspace{5mm}
\noindent \textbf{Motivating Example.}
Consider the arguably simplest instance of the data-driven optimization problem~\eqref{prob:cc general}, which estimates the worst-case value-at-risk $\sup_{\mathbb{P} \in \mathcal{F} (\theta)} \, \mathbb{P}\text{-VaR}_\varepsilon (\tilde{\xi})$ of a scalar random variable $\tilde{\xi}$ at level $\varepsilon$ from a limited set of i.i.d.~training samples $\{ \hat{\xi}_i \}_{i = 1}^N$ of $\tilde{\xi}$ under the unknown data-generating distribution $\mathbb{P}_0$ that are summarized by the empirical distribution $\hat{\mathbb{P}} = \frac{1}{N}\sum_{i = 1}^N \delta_{\hat{\xi}_i}$ at the centre of the Wasserstein ball $\mathcal{F} (\theta)$. To avoid technicalities, we assume that $\mathbb{P}_0$ is atomless. In addition, with $N_\dag = \lfloor (1-\varepsilon) N \rfloor$ and $N^\dag = \lceil (1-\varepsilon)N \rceil$ we define a distribution 
\begin{equation*}
\mathbb{P}^\dag
\; = \;
\frac{1}{N} \sum_{i = 1}^{N_\dag} \delta_{\hat{\xi}_{(i)}} + \frac{(1 - \varepsilon) N - N_\dag}{N} \delta_{\hat{\xi}_{(N^\dag)}} + \frac{N^\dag - (1 - \varepsilon) N}{N} \delta_{\hat{\xi}_{(N^\dag)} + \theta / \varepsilon} + \frac{1}{N}\sum_{i = N^\dag +1}^N \delta_{\hat{\xi}_{(i)} + \theta / \varepsilon}
\end{equation*}
to be used subsequently. Here, $\hat{\xi}_{(j)}$ denotes the $j$-th order statistic of the training samples $\{ \hat{\xi}_i \}_{i = 1}^N$.

The \emph{reliability} of the aforementioned worst-case value-at-risk, that is, the probability that it weakly exceeds the unknown true value-at-risk $\mathbb{P}_0\text{-VaR}_\varepsilon (\tilde{\xi})$, can be bounded from \emph{below} by
\begin{align*}
\mathbb{P}_0^N \left[ \sup_{\mathbb{P} \in \mathcal{F} (\theta)} \, \mathbb{P}\text{-VaR}_\varepsilon (\tilde{\xi}) \; \geq \; \mathbb{P}_0\text{-VaR}_{\varepsilon} (\tilde{\xi}) \right]
\;\; &\geq \;\;
\mathbb{P}^N_0 \left[ \mathbb{P}^\dag\text{-VaR}_\varepsilon (\tilde{\xi}) \; \geq \; \mathbb{P}_0\text{-VaR}_\varepsilon (\tilde{\xi}) \right] \\
&= \;\;
\mathbb{P}^N_0 \left[ \hat{\mathbb{P}}\text{-VaR}_\varepsilon (\tilde{\xi})  \; \geq \; \mathbb{P}_0\text{-VaR}_\varepsilon (\tilde{\xi}) - \theta / \varepsilon \right] \\
&= \;\;
1 - \mathbb{P}^N_0 \left[ \hat{\mathbb{P}}\text{-VaR}_\varepsilon (\tilde{\xi}) \; < \; \mathbb{P}_0\text{-VaR}_\varepsilon (\tilde{\xi}) - \theta / \varepsilon \right] \\
&\geq \;\;
1 - \text{exp} \left( -2N(1-\varepsilon-\mathbb{P}_0[\tilde{\xi} \leq \mathbb{P}_0\textnormal{-VaR}_{\varepsilon}(\tilde{\xi}) - \theta/\varepsilon])^2 \right),
\end{align*}
where $\mathbb{P}_0^N$ is the $N$-fold product of $\mathbb{P}_0$ that generates $\{ \hat{\xi}_i \}_{i=1}^N$. The first inequality holds since $\mathbb{P}^\dag$ is contained in $\mathcal{F} (\theta)$. The first equality holds since $\mathbb{P}^\dag\text{-VaR}_\varepsilon (\tilde{\xi}) = \hat{\mathbb{P}}\text{-VaR}_\varepsilon (\tilde{\xi}) + \theta / \varepsilon$ by construction of $\mathbb{P}^\dag$, and the last inequality is due to a standard concentration inequality for empirical quantiles (see, \emph{e.g.}, Theorem~2.3.2 of \citealt{serfling2009approximation}).

If we replace the Wasserstein ambiguity set $\mathcal{F} (\theta)$ with the ambiguity $\mathcal{G} (\theta)$ of any $\phi$-divergence, on the other hand, then we can bound the reliability from \emph{above} by
\begin{align*}
\mathbb{P}_0^N \left[ \sup_{\mathbb{P} \in \mathcal{G} (\theta)} \, \mathbb{P}\text{-VaR}_\varepsilon (\tilde{\xi}) \; \geq \; \mathbb{P}_0\text{-VaR}_{\varepsilon} (\tilde{\xi}) \right]
\;\; &= \;\;
1 - \mathbb{P}^N_0 \left[ \sup_{\mathbb{P} \in \mathcal{G} (\theta)} \, \mathbb{P}\text{-VaR}_\varepsilon (\tilde{\xi}) \; < \; \mathbb{P}_0\text{-VaR}_{\varepsilon} (\tilde{\xi}) \right] \\
\;\; &\leq \;\;
1 - \mathbb{P}^N_0 \left[ \hat{\xi}_1, \hat{\xi}_2, \ldots, \hat{\xi}_N \; < \; \mathbb{P}_0\text{-VaR}_{\varepsilon} (\tilde{\xi}) \right] \\
\;\; &\leq \;\;
1 - (1 - \varepsilon)^N.
\end{align*}
Here, the first inequality holds since all distributions in $\mathcal{G} (\theta)$ share a common support with $\hat{\mathbb{P}}$, and the second inequality follows from the definition of $\mathbb{P}_0\text{-VaR}_{\varepsilon} (\tilde{\xi})$. We highlight that this probability bound holds for every radius $\theta$ of the $\phi$-divergence ball $\mathcal{G} (\theta)$.

\begin{figure}[tb]
\begin{subfigure}{.33\textwidth}
\begin{center}
\includegraphics[width=0.95\linewidth]{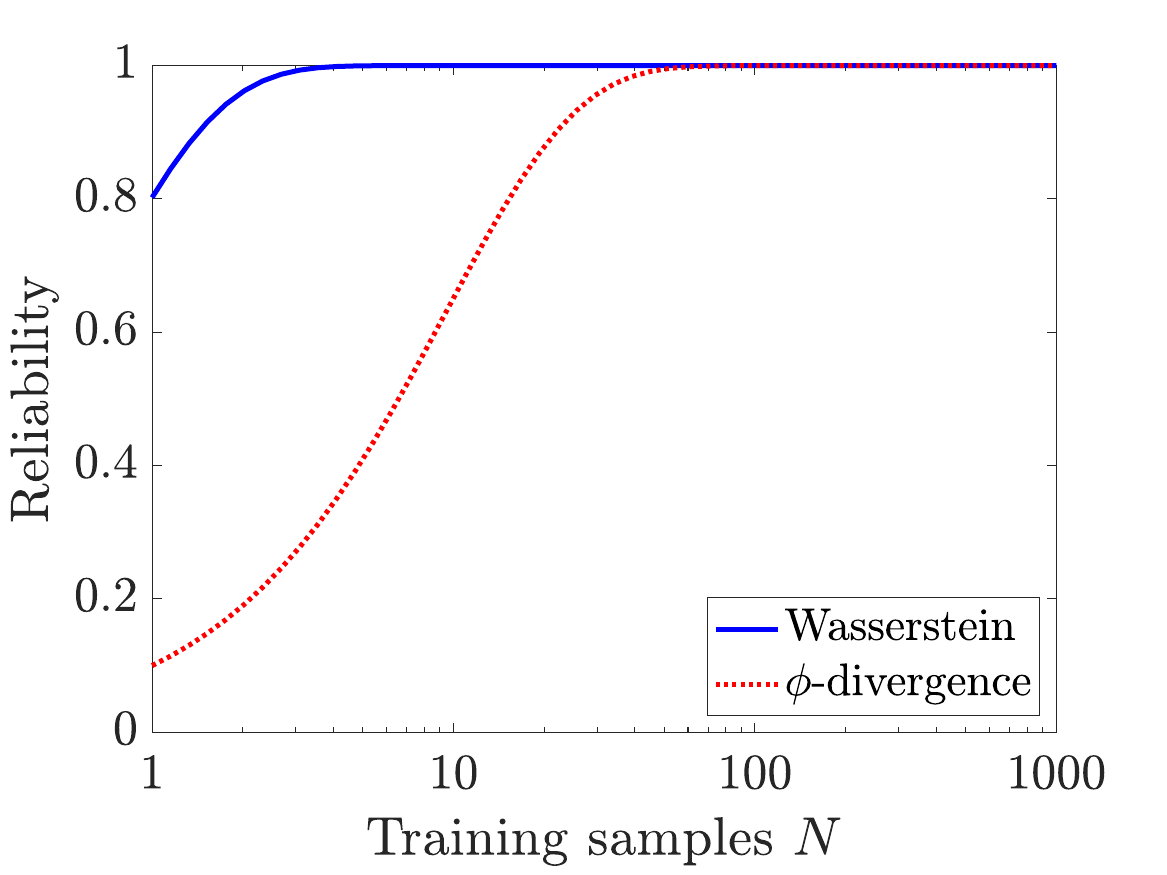}
\end{center}
\end{subfigure}%
\begin{subfigure}{.33\textwidth}
\begin{center}
\includegraphics[width=0.95\linewidth]{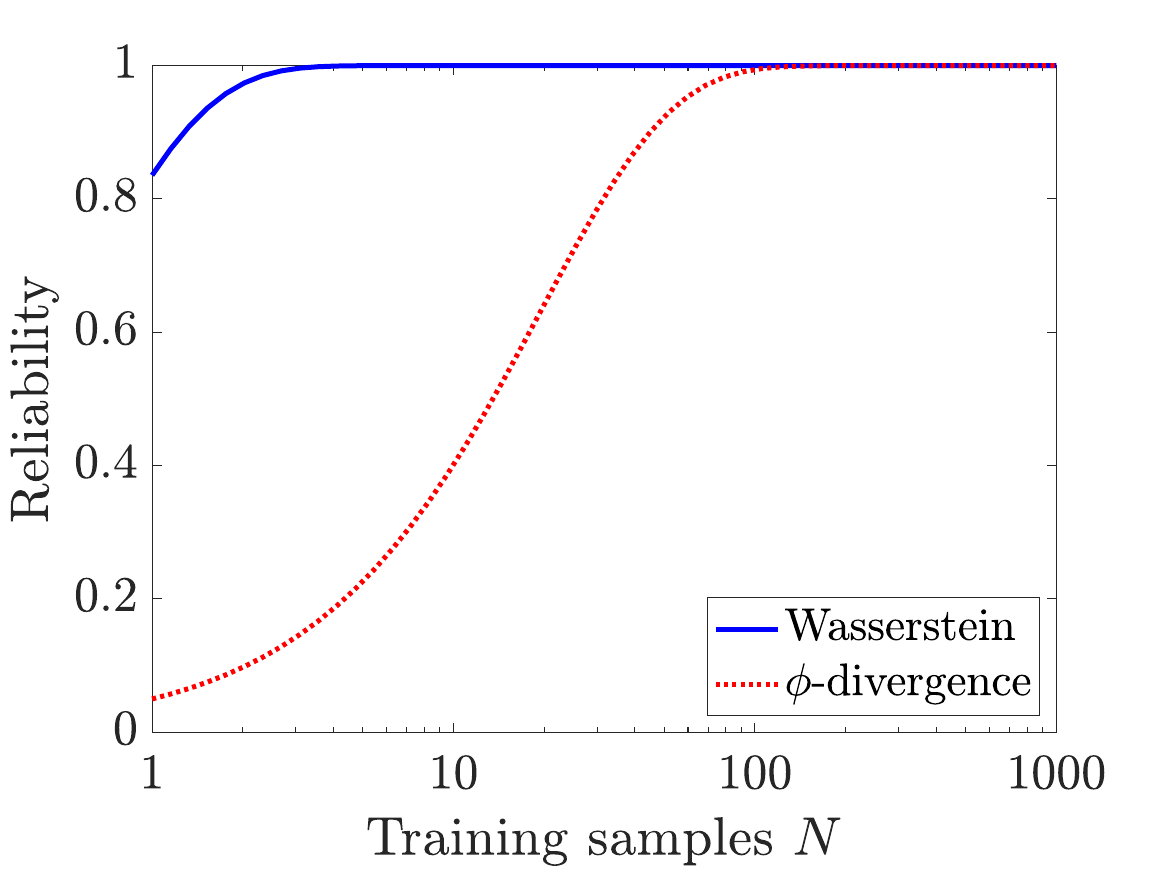}
\end{center}
\end{subfigure}	
\begin{subfigure}{.33\textwidth}
\begin{center}
\includegraphics[width=0.95\linewidth]{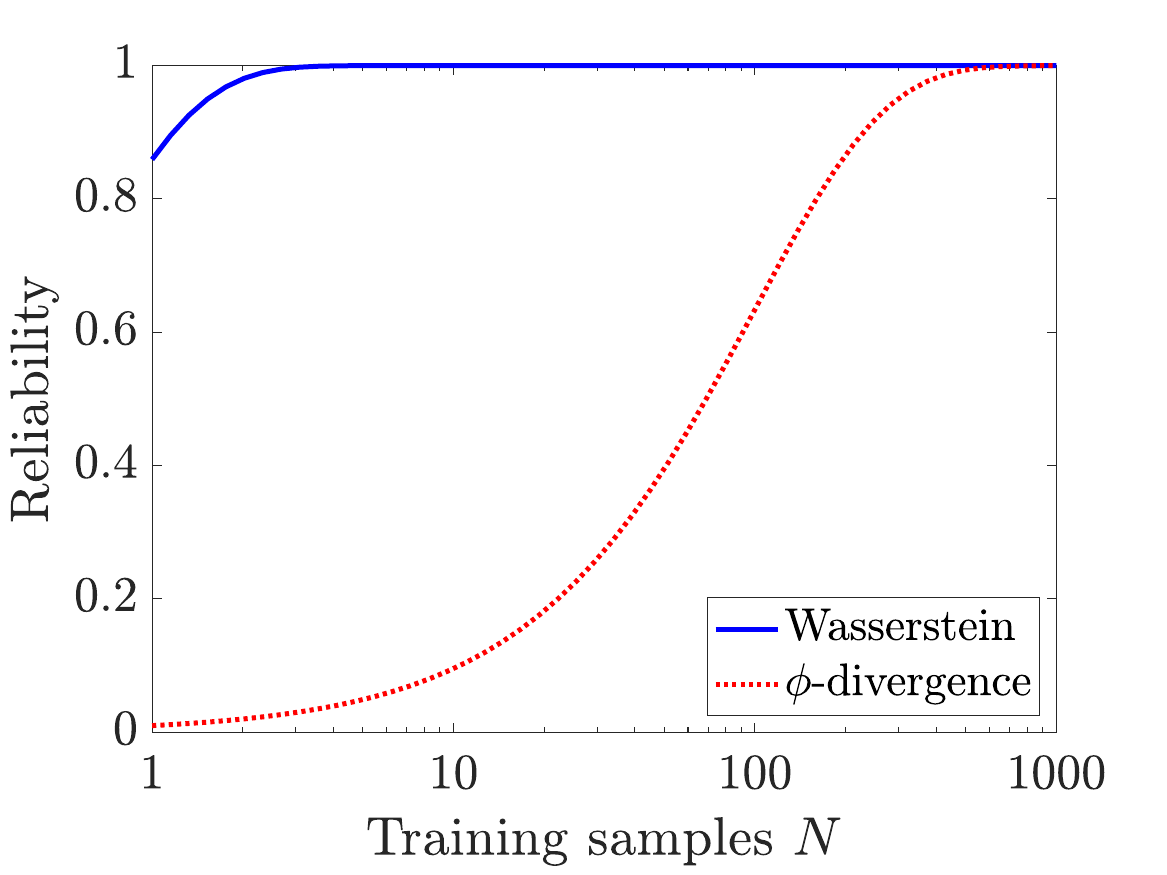}
\end{center}
\end{subfigure}
\vspace{0.2cm}
\caption{\textnormal{Reliability bounds for the Wasserstein (worst-case) and $\phi$-divergence (best-case) ambiguity sets when approximating the VaR at level $\varepsilon = 0.1$ (left), $\varepsilon =0.05$ (middle) and $\varepsilon =0.01$ (right). We choose the radius $\theta = 1/\sqrt{N}$ for the Wassestein ball (see, \emph{e.g.}, \citealt{Esfahani_Kuhn_2017}).} \label{fig:bound}}
\end{figure}

Figure~\ref{fig:bound} compares the \emph{worst-case} reliability offered by the Wasserstein ambiguity set with the \emph{best-case} reliability of the $\phi$-divergence ambiguity set for a uniform distribution over the interval~$[0, 1]$. We observe that in low-sample regimes, $\phi$-divergence ambiguity sets may underestimate the true VaR with high probability. \hfill $\clubsuit$

\vspace{5mm}
To our best knowledge, the paper of \cite{xie2020bicriteria} is the only previous work on data-driven chance constraints over Wasserstein ambiguity sets. The authors study the special class of covering problems, where the feasible region $\mathcal{X}$ satisfies $\eta \mathcal{X} \subseteq \mathcal{X}$ for every $\eta \geq 1$. This problem class encompasses, among others, portfolio optimization problems without budgetary restrictions and lot-sizing problems in the absence of setup costs. The authors prove that the resulting individual chance constrained program is NP-hard. They also demonstrate that two popular approximation schemes, the CVaR approximation as well as the scenario approximation, can perform arbitrarily poorly for classical individual chance constraints, that is, when the Wasserstein radius is $\theta = 0$. Based on this insight, the authors propose a bicriteria approximation scheme for covering problems with classical as well as distributionally robust individual chance constraints over moment and Wasserstein ambiguity sets. This bicriteria approximation scheme determines solutions that trade off a higher risk threshold $\varepsilon' > \varepsilon$ in the chance constraint with a smaller optimality gap $\varepsilon' / (\varepsilon' - \varepsilon)$. This is achieved by solving a tractable convex relaxation of the chance constrained problem (using, \emph{e.g.}, a Markovian or Bernstein generator) and subsequently scaling the solution to this relaxation so that it becomes feasible for the chance constraint with the higher risk threshold $\varepsilon'$. By design, the performance guarantee of the bicriteria approximation scheme becomes weaker (and eventually trivial) as the selected risk threshold $\varepsilon'$ \mbox{approaches the risk threshold $\varepsilon$ of the original problem formulation.}

In this paper, we study distributionally robust chance constrained programs over the Wasserstein ambiguity set~\eqref{set:Wasserstein}. We derive deterministic reformulations for individual chance constrained programs, where $\mathcal{S}(\bm{x}) = \{ \bm{\xi} \in \mathbb{R}^K \mid \bm{a} (\bm{\xi})^\top \bm{x} < b (\bm{\xi}) \}$ for affine functions $\bm{a}(\cdot) : \mathbb{R}^K \rightarrow \mathbb{R}^L$ and $b(\cdot) : \mathbb{R}^K \rightarrow \mathbb{R}$, as well as for joint chance constrained programs with right-hand side uncertainty, where $\mathcal{S}(\bm{x}) = \{\bm{\xi} \in \mathbb{R}^K \mid \bm{A} \bm{x} < \bm{b} (\bm{\xi}) \}$ for $\bm{A} \in \mathbb{R}^{M \times L}$ and an affine function $\bm{b} : \mathbb{R}^K \rightarrow \mathbb{R}^M$. Our reformulations are mixed-integer conic programs that reduce to mixed-integer linear programs when the norm $\left \lVert \cdot \right \rVert$ on $\mathbb{R}^K$ is the $1$-norm or the $\infty$-norm. 

While preparing this paper for publication, we became aware of the independent work by \cite{xie2019distributionally}, which derives similar reformulations for distributionally individual and joint chance constraints over Wasserstein ambiguity sets. In contrast to our work, however, \cite{xie2019distributionally} assumes that each safety condition $\bm{a}_m^\top \bm{x} < b_m (\bm{\xi})$, $m \in [M]$, in the joint chance constraint depends on a subvector of $\bm{\xi}$, and that these subvectors are pairwise disjoint for different safety conditions. In other words, different safety conditions of the joint chance constraints studied by \cite{xie2019distributionally} must depend on different random variables. Furthermore, the reformulations of \cite{xie2019distributionally} are derived via duality theory, whereas our reformulations directly leverage the structural insights into the worst-case distributions. This enables us to keep our reformulations largely independent of the selected ground metric for the Wasserstein ball, which opens up possibilities to incorporate other cost functions in our definition of the Wasserstein distance. 
Since the initial submission of this paper, our exact reformulation for data-driven chance constrained program over Wasserstein balls has been further studied and tightened; see, for instance, \cite{ho2020strong, ho2021distributionally}, \cite{shen2021convex} and \cite{zhang2021building}. Along with these theoretical extensions, our reformulation has also been applied in several domains, including risk sharing in finance \citep{chen2021sharing}, network design for humanitarian operations \citep{jiang2021distributionally} and optimal power flows in energy systems \citep{arrigo2022wasserstein}.

\vspace{5mm}
\noindent \textbf{Notation.}
We use boldface uppercase and lowercase letters to denote matrices and vectors, respectively. Special vectors of appropriate dimensions include $\bm{0}$ and $ \bm{e} $, which respectively correspond to the zero vector and the vector of all ones. We denote by $\|\cdot\|_*$ the dual norm of a general norm $\|\cdot\|$. We use the shorthand $ [N] = \left\{1,2,\ldots,N\right\} $ to represent the set of all integers up to $ N $. Given a (possibly fractional) real number $\ell\in [0,N]$, we define the partial sum of the $\ell$ first values in $\{k_i\}_{i \in [N]}$ as $\sum_{i = 1}^{\ell} k_i = \sum_{i = 1}^{\lfloor \ell \rfloor} k_i + (\ell - \lfloor \ell \rfloor) k_{\lfloor \ell \rfloor + 1}$.  Random vectors are denoted by tilde signs ({\em e.g.}, $\bmt{\xi}$), while their realizations are denoted by the same symbols without tildes ({\em e.g.}, $\bm{\xi}$). Given a random vector $\bmt{\xi}$ governed by a distribution $\mathbb{P}$, a measurable loss function $\ell (\bm{\xi})$ and a risk threshold $\varepsilon \in (0, 1)$, the value-at-risk (VaR) of $\ell (\bm{\xi})$ at level $\varepsilon$ is defined as $\mathbb{P}\text{-VaR}_{\varepsilon} (\ell (\bm{\xi})) = \inf\{\gamma \in \mathbb{R} \mid \mathbb{P}[\gamma \leq \ell(\bmt{\xi})] \leq \varepsilon\}$.

\section{Exact Reformulation of Data-Driven Chance Constraints}\label{sec:exact_reformulation}

Section~\ref{sec:uq_wasserstein} reviews a previously established result on the quantification of uncertainty over Wasserstein balls. We use this result to derive an exact reformulation of generic data-driven chance constrained programs in Section~\ref{sec:ref_generic}. We finally specialize this generic reformulation to the subclasses of data-driven individual chance constrained programs as well as data-driven joint chance constrained programs with right-hand side uncertainty in Sections~\ref{sec:ref_indiv_cc} and~\ref{sec:ref_joint_cc}, respectively.

\subsection{Uncertainty Quantification over Wasserstein Balls}\label{sec:uq_wasserstein}

Consider an open safety set $\mathcal{S} \subseteq \mathbb{R}^K $, and denote by $\bar{\mathcal{S}} = \mathbb{R}^K \setminus \mathcal{S}$ its closed complement. The uncertainty quantification problem 
\begin{equation}
\label{prob:uncertainty quantification}
\sup_{\mathbb{P} \in \mathcal{F}(\theta)} \mathbb{P}[\bmt{\xi} \notin \mathcal{S}]
\end{equation}
computes the worst (largest) probability of the system under consideration being unsafe, which is the case whenever the random vector $\bmt{\xi}$ attains a value in the unsafe set $\bar{\mathcal{S}}$. Throughout the rest of the paper, we exclude trivial special cases and assume that $\theta > 0$ and $\varepsilon \in (0, 1)$.

To solve the uncertainty quantification problem~\eqref{prob:uncertainty quantification}, denote by $\mathbf{dist}(\bmh{\xi}_i, \bar{\mathcal{S}})$ the distance of the $i^\text{th}$ data point $\bmh{\xi}_i \in \mathbb{R}^K$ of the empirical distribution $\hat{\mathbb{P}}$ to the unsafe set $\bar{\mathcal{S}}$. This distance is based on a norm $\left \lVert \cdot \right \rVert$, which we keep generic at this stage. Without loss of generality, we assume that the data points $\{\bmh{\xi}_i\}_{i \in [N]} $ are ordered in increasing distance to $\bar{\mathcal{S}}$, that is, $\mathbf{dist}(\bmh{\xi}_i, \bar{\mathcal{S}}) \leq \mathbf{dist}(\bmh{\xi}_j, \bar{\mathcal{S}})$ for all $1 \leq i \leq j \leq N$. We also assume that $\mathbf{dist}(\bmh{\xi}_i, \bar{\mathcal{S}}) = 0$ (that is, the data point $\bmh{\xi}_i$ is unsafe) if and only if $i \in [I]$, where $I = 0$ if $\mathbf{dist}(\bmh{\xi}_i, \bar{\mathcal{S}}) > 0$ for all $i \in [N]$. Finally, we denote by $\bm{\xi}^\star_i \in \bar{\mathcal{S}}$ an unsafe point that is closest to the data point $\bmh{\xi}_i$, $i \in [N]$, in terms of the distance $\mathbf{dist} (\bmh{\xi}_i, \bar{\mathcal{S}})$.

\cite{blanchet2019quantifying} as well as \cite{Gao_Kleywegt_2016} have characterized the solution to the uncertainty quantification problem~\eqref{prob:uncertainty quantification} in closed form. To keep our paper self-contained, we reproduce their findings without proof in Theorem~\ref{thm:uncertainty-quantification} below.

\begin{theorem}\label{thm:uncertainty-quantification}
Let $j^\star = \max \, \{j \in [N] \cup \{ 0 \} \mid \sum_{i = 1}^j \mathbf{dist}(\bmh{\xi}_i, \bar{\mathcal{S}}) \leq \theta N \}$. The uncertainty quantification problem~\eqref{prob:uncertainty quantification} is solved by a worst-case distribution $\mathbb{P}^\star \in \mathcal{F} (\theta)$ that is characterized as follows:
\begin{enumerate}
\item[(i)] If $j^\star = N$, then $\sup\limits_{\mathbb{P} \in \mathcal{F}(\theta)} \mathbb{P}[\bmt{\xi} \notin \mathcal{S}] \; = \; \mathbb{P}^\star [\bmt{\xi} \notin \mathcal{S}] \; = \; 1$ for
\begin{equation*}
\mathbb{P}^\star \; = \; \dfrac{1}{N} \sum_{i = 1}^I \delta_{\bmh{\xi}_i} \; + \; \dfrac{1}{N} \sum_{i = I+1}^N \delta_{\bm{\xi}^\star_i}.
\end{equation*}
\item[(ii)] If $j^\star < N$, then $\sup\limits_{\mathbb{P} \in \mathcal{F}(\theta)} \mathbb{P}[\bmt{\xi} \notin \mathcal{S}] \; = \; \mathbb{P}^\star [\bmt{\xi} \notin \mathcal{S}] \; = \; (j^\star+ p^\star)/N$ for
\begin{equation*}
\mathbb{P}^\star \; = \; \dfrac{1}{N} \sum_{i = 1}^I \delta_{\bmh{\xi}_i} \; + \; \dfrac{1}{N} \sum_{i = I+1}^{j^\star} \delta_{\bm{\xi}^\star_i} \; + \; \dfrac{p^\star}{N} \delta_{\bm{\xi}^\star_{j^\star+1}} \; + \;  \dfrac{1-p^\star}{N}\delta_{\bmh{\xi}_{j^\star+1}} \; + \; \dfrac{1}{N} \sum_{i = j^\star+2}^N  \delta_{\bmh{\xi}_i},
\end{equation*}
where $p^\star = (\theta N - \sum_{i = 1}^{j^\star} \mathbf{dist}(\bmh{\xi}_i, \bar{\mathcal{S}})) / \mathbf{dist}(\bmh{\xi}_{j^\star+1}, \bar{\mathcal{S}})$.
\end{enumerate}
\end{theorem}

\begin{figure}[tb]
\begin{subfigure}{.5\textwidth}
\begin{center}
\includegraphics[width=0.75\linewidth]{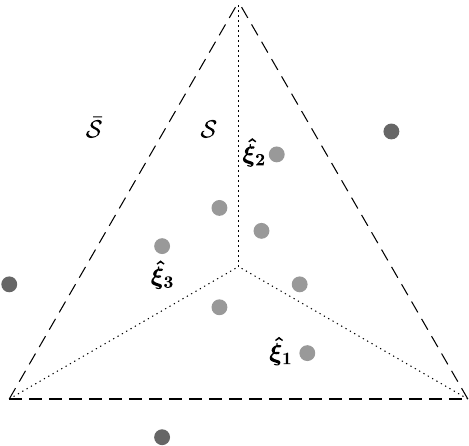}
\end{center}
\end{subfigure}%
\begin{subfigure}{.5\textwidth}
\begin{center}
\includegraphics[width=0.75\linewidth]{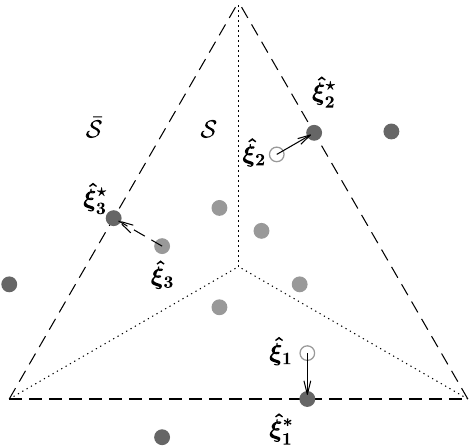}
\end{center}
\end{subfigure}
\vspace{0.2cm}
	
\caption{{\textnormal{Empirical and worst-case distributions. The left graph visualizes the empirical distribution $\hat{\mathbb{P}}$, whose light grey (dark grey) data points are contained in (outside of) the safety set $\mathcal{S}$ shown as an equilateral triangle (dashed lines). The right graph shows the corresponding worst-case distribution $\mathbb{P}^\star$, which moves the data points $\bmh{\xi}_1$ and $\bmh{\xi}_2$ entirely as well as the data point $\bmh{\xi}_3$ partially to the unsafe set $\bar{\mathcal{S}}$. Each transported data point is projected onto the boundary of the closest halfspace defining the safety set $\mathcal{S}$.}} \label{fig:greedy}}
\end{figure}

Intuitively speaking, the worst-case distribution $\mathbb{P}^\star$ in Theorem~\ref{thm:uncertainty-quantification} transports the training dataset $\{\bmh{\xi}_i\}_{i \in [N]}$ to the unsafe set $\bar{\mathcal{S}}$ in a greedy fashion, see Figure~\ref{fig:greedy}. The data points $\bmh{\xi}_1,\dots,\bmh{\xi}_I$ are already unsafe and hence do not need to be transported. The subsequent data points $\bmh{\xi}_{I+1}, \ldots, \bmh{\xi}_{j^\star + 1}$ are closest to the unsafe set and are thus transported from $\mathcal{S}$ to $\bar{\mathcal{S}}$. Due to the limited transportation budget $\theta$, the data point $\bmh{\xi}_{j^\star + 1}$ is only partially transported. The safe samples $\bmh{\xi}_{j^\star + 2}, \ldots \bmh{\xi}_N$, finally, are too far away from the unsafe set $\bar{\mathcal{S}}$ and are thus left unchanged. Note that the distribution characterized in Theorem~\ref{thm:uncertainty-quantification} may not be the only distribution that solves problem~\eqref{prob:uncertainty quantification}.

\subsection{Reformulation of Generic Chance Constraints}\label{sec:ref_generic}

We now develop deterministic reformulations for the distributionally robust chance constrained program~\eqref{prob:cc general}. To this end, we focus on the ambiguous chance constraint
\begin{equation}
\label{prob:worst-case cc}
\sup_{\mathbb{P} \in \mathcal{F}(\theta)} \mathbb{P}[\bmt{\xi} \notin \mathcal{S}(\bm{x})] \leq \varepsilon.
\end{equation}
For any fixed decision $\bm x\in \mathcal{X}$, we let $\mathcal{S}(\bm{x})$ be an arbitrary open safety set, and we denote by $\bar{\mathcal{S}}(\bm{x})$ its closed complement, which comprises all unsafe scenarios. 
Every fixed training dataset $\{\bmh{\xi}_i\}_{i \in [N]}$ then induces a (decision-dependent) permutation $\bm{\pi}(\bm{x})$ of $[N]$ that orders the training samples in increasing distance to the unsafe set, that is, 
\begin{equation*}
\mathbf{dist}(\bmh{\xi}_{\pi_1(\bm{x})}, \bar{\mathcal{S}}(\bm{x})) \; \leq \; \mathbf{dist}(\bmh{\xi}_{\pi_2(\bm{x})}, \bar{\mathcal{S}}(\bm{x})) 
\; \leq \; \cdots \; \leq \; \mathbf{dist}(\bmh{\xi}_{\pi_N(\bm{x})}, \bar{\mathcal{S}}(\bm{x})).
\end{equation*}

We first show that a fixed decision $\bm{x}$ satisfies the ambiguous chance constraint~\eqref{prob:worst-case cc} over the Wasserstein ambiguity set~\eqref{set:Wasserstein} if and only if the partial sum of the $\varepsilon N$ smallest transportation distances to the unsafe set multiplied by the mass $1/N$ of a training sample exceeds~$\theta$.

\begin{theorem}\label{thm:cc equivalent}
For any fixed decision $\bm{x}\in \mathcal{X}$, the ambiguous chance constraint~\eqref{prob:worst-case cc} over the Wasserstein ambiguity set~\eqref{set:Wasserstein} is equivalent to the deterministic inequality
\begin{equation}
\label{equivalence:theta positive}
\dfrac{1}{N}\sum_{i = 1}^{\varepsilon N} \mathbf{dist}(\bmh{\xi}_{\pi_i(\bm{x})}, \bar{\mathcal{S}}(\bm{x})) \ge \theta.
\end{equation}
\end{theorem} 

The left-hand side of~\eqref{equivalence:theta positive} can be interpreted as the minimum cost of moving a fraction $\varepsilon$ of the training samples to the unsafe set. If this cost exceeds the prescribed transportation budget $\theta$, then no distribution in the Wasserstein ambiguity set can assign the unsafe set a probability of more than $\varepsilon$, which means that the distributionally robust chance constraint~\eqref{prob:worst-case cc} is satisfied.

\noindent \emph{Proof of Theorem~\ref{thm:cc equivalent}.} $\;$
From Theorem~\ref{thm:uncertainty-quantification} we know that the worst-case distribution $\mathbb{P}^\star$ is an optimal solution (not necessarily unique) to the maximization problem embedded in the left-hand side of the ambiguous chance constraint~\eqref{prob:worst-case cc}. We thus conclude that the constraint~\eqref{prob:worst-case cc} is satisfied if and only if $\mathbb{P}^\star [\bmt{\xi} \notin \mathcal{S}(\bm{x})] \leq \varepsilon$ for  $\mathbb{P}^\star$ defined in the statement of that theorem.

In case (\emph{i}) of Theorem~\ref{thm:uncertainty-quantification}, the ambiguous chance constraint~\eqref{prob:worst-case cc} is violated since $\mathbb{P}^\star [\bmt{\xi} \notin \mathcal{S}(\bm{x})] = 1$ while $\varepsilon < 1$ by assumption. At the same time, since $j^\star = N$, we have $\frac{1}{N} \sum_{i = 1}^{N} \mathbf{dist}(\bmh{\xi}_{\pi_i(\bm{x})}, \bar{\mathcal{S}}(\bm{x})) \leq \theta$. If this inequality is strict, then~\eqref{equivalence:theta positive} is violated as desired since $\frac{1}{N} \sum_{i = 1}^{\varepsilon N} \mathbf{dist}(\bmh{\xi}_{\pi_i(\bm{x})}, \bar{\mathcal{S}}(\bm{x})) \leq \frac{1}{N} \sum_{i = 1}^{N} \mathbf{dist}(\bmh{\xi}_{\pi_i(\bm{x})}, \bar{\mathcal{S}}(\bm{x}))$. If the inequality is satisfied as an equality, on the other hand, we know that $\mathbf{dist}(\bmh{\xi}_{\pi_N(\bm{x})}, \bar{\mathcal{S}}(\bm{x})) > 0$ since $\theta > 0$ by assumption and $\mathbf{dist}(\bmh{\xi}_{\pi_i(\bm{x})}, \bar{\mathcal{S}}(\bm{x})) \leq \mathbf{dist}(\bmh{\xi}_{\pi_j(\bm{x})}, \bar{\mathcal{S}}(\bm{x}))$ for all $i \leq j$ by construction of the re-ordering $\bm{\pi}(\bm{x})$. Thus, since $\varepsilon < 1$ by assumption, we have $\frac{1}{N} \sum_{i = 1}^{\varepsilon N} \mathbf{dist}(\bmh{\xi}_{\pi_i(\bm{x})}, \bar{\mathcal{S}}(\bm{x})) < \frac{1}{N} \sum_{i = 1}^{N} \mathbf{dist}(\bmh{\xi}_{\pi_i(\bm{x})}, \bar{\mathcal{S}}(\bm{x})) = \theta$ and equation~\eqref{equivalence:theta positive} is violated as desired.

In case (\emph{ii}) of Theorem~\ref{thm:uncertainty-quantification}, we have $\mathbb{P}^\star [\bmt{\xi} \notin \mathcal{S} (\bm{x})] = (j^\star+ p^\star)/N$ with $j^\star = \max \, \{j \in [N - 1] \cup \{0\} \mid \sum_{i = 1}^j \mathbf{dist}(\bmh{\xi}_{\pi_i (\bm{x})}, \bar{\mathcal{S}} (\bm{x})) \leq \theta N \}$ as well as $p^\star = (\theta N - \sum_{i = 1}^{j^\star} \mathbf{dist}(\bmh{\xi}_{\pi_i (\bm{x})}, \bar{\mathcal{S}} (\bm{x}))) / \mathbf{dist}(\bmh{\xi}_{\pi_{j^\star + 1} (\bm{x})}, \bar{\mathcal{S}} (\bm{x}))$. We claim that $j^\star+ p^\star$ is the optimal value of the bivariate mixed-integer optimization problem
\begin{equation}\label{eq:bivariate_mixed_integer}
\begin{array}{cll}
\displaystyle \max_{j, p} & \displaystyle j + p \\
{\rm s.t.} & \displaystyle  \sum_{i = 1}^{j} \mathbf{dist}(\bmh{\xi}_{\pi_i (\bm{x})}, \bar{\mathcal{S}} (\bm{x})) + p \cdot \mathbf{dist}(\bmh{\xi}_{\pi_{j + 1} (\bm{x})}, \bar{\mathcal{S}} (\bm{x})) \leq \theta N \\[4mm]
& \displaystyle j \in [N - 1] \cup \{0\},~0 \leq p < 1.
\end{array}
\end{equation}
Indeed, the solution $(j, p) = (j^\star, p^\star)$ is feasible in~\eqref{eq:bivariate_mixed_integer} by definition of $j^\star$ and $p^\star$. Moreover, we have $j + p < j^\star + p^\star$ for any other feasible solution $(j, p)$ that satisfies $j = j^\star$ and $p \neq p^\star$. Assume now that the optimal solution $(j, p)$ to~\eqref{eq:bivariate_mixed_integer} would satisfy $j > j^\star$. Any such solution would violate the first constraint since $\sum_{i = 1}^j \mathbf{dist}(\bmh{\xi}_{\pi_i (\bm{x})}, \bar{\mathcal{S}} (\bm{x})) > \theta N$ by definition of $j^\star$ while $p \geq 0$. Similarly, any solution $(j, p)$ with $j < j^\star$ cannot be optimal in~\eqref{eq:bivariate_mixed_integer} since $j \leq j^\star - 1$ while $p < p^\star + 1$.

\begin{figure}[tb]
\begin{subfigure}{.5\textwidth}
\begin{center}
\includegraphics[width=0.9\linewidth]{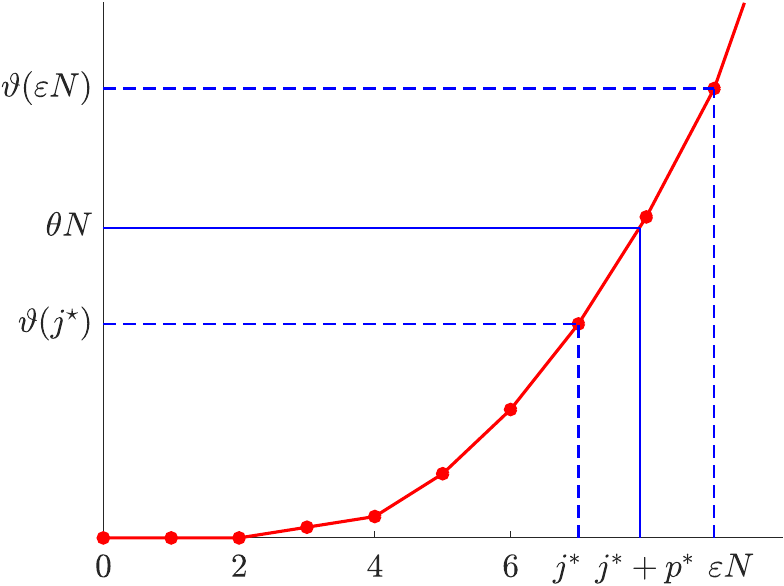}
\end{center}
\end{subfigure}%
\begin{subfigure}{.5\textwidth}
\begin{center}
\includegraphics[width=0.9\linewidth]{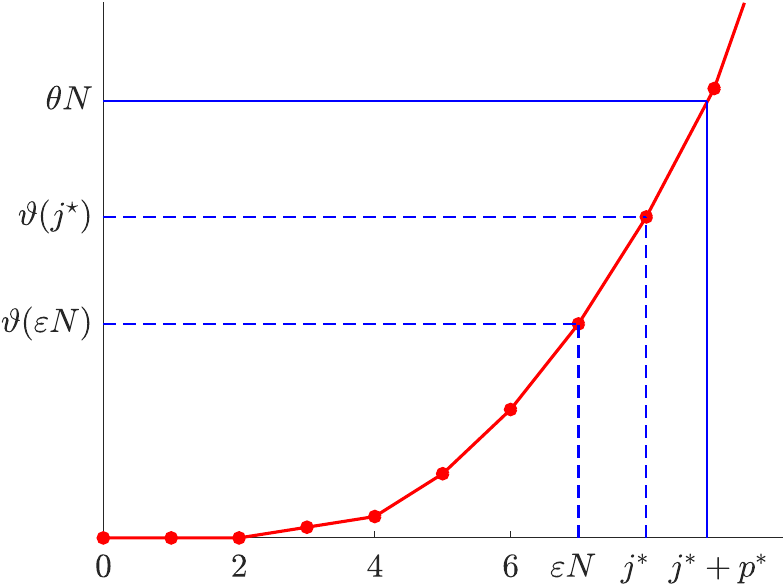}
\end{center}
\end{subfigure}
\vspace{0.2cm}
\caption{{\textnormal{Relationship between $\vartheta$ and $\vartheta^{-1}$. The left graph shows a feasible solution $\bm{x}$ satisfying the ambiguous chance constraint~\eqref{prob:worst-case cc}; in this case, we have $\vartheta (\varepsilon N) \geq \theta N$. The infeasible solution $\bm{x}'$ in the right graph, on the other hand, violates the ambiguous chance constraint~\eqref{prob:worst-case cc}, and we have $\vartheta (\varepsilon N) < \theta N$.}} \label{fig:inverse_function}}
\end{figure}

We can re-express problem~\eqref{eq:bivariate_mixed_integer} as the univariate discrete optimization problem
\begin{equation*}
\max \bigg\{ j \in [0, N] ~\bigg|~ \sum_{i = 1}^{\lfloor j \rfloor} \mathbf{dist}(\bmh{\xi}_{\pi_i (\bm{x})}, \bar{\mathcal{S}} (\bm{x})) \; + \; (j - \lfloor j \rfloor) \cdot \mathbf{dist}(\bmh{\xi}_{\pi_{\lfloor j \rfloor + 1} (\bm{x})}, \bar{\mathcal{S}} (\bm{x})) \leq \theta N \bigg\}.
\end{equation*}
Using our definition of partial sums, we observe that this problem is equivalent to
\begin{equation*}
\max \bigg\{ j \in [0, N] ~\bigg|~ \sum_{i = 1}^j \mathbf{dist}(\bmh{\xi}_{\pi_i (\bm{x})}, \bar{\mathcal{S}} (\bm{x})) \leq \theta N \bigg\}.
\end{equation*}
By construction, the mapping $\vartheta (j) = \sum_{i = 1}^j \mathbf{dist}(\bmh{\xi}_{\pi_i (\bm{x})}, \bar{\mathcal{S}} (\bm{x}))$, $j \in [0, N]$, is continuous and monotonically nondecreasing. It therefore affords the right inverse $\vartheta^{-1} (t) = \max \{ j \in [0, N] \mid \vartheta(j) \leq t \}$ that satisfies $\vartheta \circ \vartheta^{-1} (t) = t$ for all $t \in [0, \vartheta(N)]$. Figure~\ref{fig:inverse_function} visualizes the relationship between $\vartheta$ and $\vartheta^{-1}$. We thus conclude that the ambiguous chance constraint~\eqref{prob:worst-case cc} is satisfied if and only if
\begin{align*}
\max \bigg\{ j \in [0, N] ~\bigg|~ \sum_{i = 1}^j \mathbf{dist}(\bmh{\xi}_{\pi_i (\bm{x})}, \bar{\mathcal{S}} (\bm{x})) \leq \theta N \bigg\} \leq \varepsilon N \quad
&\Longleftrightarrow \quad \max \{ j \in [0, N] ~|~ \vartheta (j) \leq \theta N \} \leq \varepsilon N \\
&\Longleftrightarrow \quad \vartheta^{-1} (\theta N) \leq \varepsilon N \\
&\Longleftrightarrow \quad \theta N \leq \vartheta (\varepsilon N),
\end{align*}
where the last equivalence follows from $\vartheta \circ \vartheta^{-1} (\theta N) = \theta N$, which holds because $\theta N \leq \vartheta(N)$ for $j^\star < N$, as well as the fact that $\vartheta$ is monotonically nondecreasing. By definition, the right-hand side of the last equivalence holds if and only if~\eqref{equivalence:theta positive} in the statement of the theorem is satisfied.
\hfill \Halmos
\endproof

\begin{remark}
We emphasize that the inequality~\eqref{equivalence:theta positive} fails to be equivalent to the ambiguous chance constraint~\eqref{prob:worst-case cc} when $\theta = 0$, in which case the Wasserstein ball collapses to the singleton set $\mathcal{F}(0) = \{\hat{\mathbb{P}}\}$. To see this, suppose that $\bmh{\xi}_{\pi_i(\bm{x})} \in \bar{\mathcal{S}}(\bm{x})$ for all $i=1,\ldots,I$ and $\bmh{\xi}_{\pi_i(\bm{x})} \in \mathcal{S}(\bm{x})$ for all $i=I+1,\ldots,N$, where $I\ge 1$. If $\varepsilon < I/N$, then the chance constraint~\eqref{prob:worst-case cc} is violated because
\begin{equation*}
\hat{\mathbb{P}}[\bmt{\xi} \notin \mathcal{S}(\bm{x})] = \frac{I}{N}>\varepsilon,
\end{equation*}
while the inequality~\eqref{equivalence:theta positive} holds trivially because $\sum_{i = 1}^{\varepsilon N} \mathbf{dist}(\bmh{\xi}_{\pi_i(\bm{x})}, \bar{\mathcal{S}}(\bm{x})) \ge 0$.
\end{remark}

Theorem~\ref{thm:cc equivalent} establishes that a decision $\bm{x} \in \mathcal{X}$ satisfies the ambiguous chance constraint~\eqref{prob:worst-case cc} if and only if 
the sum of the $\varepsilon N$ smallest distances of the training samples to the unsafe set $\bar{\mathcal{S}}(\bm{x})$ weakly exceeds $\theta N$. This result is of computational interest because the sum of the $\varepsilon N$ smallest out of $N$ real numbers is concave in those real numbers (while being convex in $\varepsilon$). This reveals that the  constraint~\eqref{equivalence:theta positive} is convex in the decision-dependent distances $\{\mathbf{dist}(\bmh{\xi}_{i}, \bar{\mathcal{S}}(\bm{x}))\}_{i \in [N]}$. In the remainder we develop an efficient reformulation of this convex constraint that does not require an enumeration of all possible sums of $\varepsilon N$ different distances between the training samples and the unsafe set. This reformulation is based on the following auxiliary lemma.

\begin{lemma}\label{lem:sum of smallest}
For any $\varepsilon \in (0,1)$, the sum of the $\varepsilon N$ smallest out of $N$ real numbers $k_1,\dots,k_N$ coincides with the optimal value of the linear program
$$
\begin{array}{cll}
\displaystyle \max_{\bm{s}, t} & \varepsilon N t - \bm{e}^\top\bm{s} \\
{\rm s.t.} & k_i \geq t - s_i &~\forall i \in [N] \\
& \bm{s} \geq \bm{0}.
\end{array}
$$
\end{lemma}
\noindent \emph{Proof of Lemma~\ref{lem:sum of smallest}.} $\;$
By definition, the sum of the $\varepsilon N$ smallest elements of the set $\{k_1,\dots,k_N\}$ corresponds to the optimal value of the (manifestly feasible) linear program
\begin{equation*}
\begin{array}{cl}
\displaystyle \min_{\bm{v}} & \displaystyle \sum_{i \in [N]} k_i v_i \\ \text{s.t.} & \bm{0} \leq \bm{v} \leq \bm{e}, ~\bm{e}^\top\bm{v} = \varepsilon N.
\end{array}
\end{equation*}
The claim now follows from strong linear programming duality.
\hfill \Halmos
\endproof 

Armed with Theorem~\ref{thm:cc equivalent} and Lemma~\ref{lem:sum of smallest}, we are now ready to reformulate the chance constrained program~\eqref{prob:cc general} as a deterministic optimization problem.

\begin{theorem}\label{thm:cc-deterministic}
The chance constrained program~\eqref{prob:cc general} is equivalent to 
\begin{equation}\label{prob:cc reformulation}
\begin{array}{cll}
\displaystyle \min_{\bm{s}, t, \bm{x}} & \bm{c}^\top\bm{x} \\
{\rm s.t.} & \varepsilon N t - \bm{e}^\top\bm{s} \geq \theta N \\
& \mathbf{dist}(\bmh{\xi}_i, \bar{\mathcal{S}}(\bm{x})) \geq t - s_i &~\forall i \in [N] \\
& \bm{s} \geq \bm{0}, ~\bm{x} \in \mathcal{X}.
\end{array}
\end{equation}
\end{theorem}
\noindent \emph{Proof of Theorem~\ref{thm:cc-deterministic}.} $\;$
The claim follows immediately by using Theorem~\ref{thm:cc equivalent} to reformulate the chance constraint~\eqref{prob:worst-case cc} as the inequality~\eqref{equivalence:theta positive}, using Lemma~\ref{lem:sum of smallest} to express the left-hand side of~\eqref{equivalence:theta positive} as a linear maximization problem and substituting the resulting constraint back into~\eqref{prob:cc general}.
\hfill \Halmos
\endproof 

We emphasize that the reformulation offered by Theorem~\ref{thm:cc-deterministic} is independent of the selected ground metric $\mathbf{dist} (\cdot, \cdot)$. In the remainder, we assume that the ground metric is based on a norm $\lVert \cdot \rVert$.

\subsection{Reformulation of Individual Chance Constraints}\label{sec:ref_indiv_cc}

Assume now that problem~\eqref{prob:cc general} accommodates an individual chance constraint defined through the safety set $\mathcal{S}(\bm{x}) = \{\bm{\xi} \in \mathbb{R}^K \mid (\bm{A}\bm{\xi} + \bm{a})^\top \bm{x} < \bm{b}^\top\bm{\xi} + b\}$. Individual chance constrained programs have been studied, among others, in network design \citep{wang2007beta}, vehicle routing \citep{gounaris2013robust, ghosal2020distributionally} and portfolio optimization \citep{rujeerapaiboon2016robust, dert2000optimal}. By Lemma~\ref{lem:distance to the union of closed half-spaces} in the appendix, we have
$$\mathbf{dist}(\bmh{\xi}_i, \bar{\mathcal{S}}(\bm x)) = \dfrac{((\bm{b} - \bm{A}^\top\bm{x})^\top\bmh{\xi}_i + b - \bm{a}^\top\bm{x})^+}{\|\bm{b} - \bm{A}^\top\bm{x}\|_*} ~~\forall i \in [N],
$$
where we adopt the convention that $0 / 0 = 0$, and thus Theorem~\ref{thm:cc-deterministic} allows us to reformulate problem~\eqref{prob:cc reformulation} as the deterministic optimization problem
\begin{equation}\label{prob:individual cc reformulation}
\begin{array}{cll}
\displaystyle \min_{\bm{s}, t, \bm{x}} & \bm{c}^\top\bm{x} \\
{\rm s.t.} & \varepsilon N t - \bm{e}^\top\bm{s} \geq \theta N \\
& \dfrac{((\bm{b} - \bm{A}^\top\bm{x})^\top\bmh{\xi}_i + b - \bm{a}^\top\bm{x})^+}{\|\bm{b} - \bm{A}^\top\bm{x}\|_*} \geq t - s_i &~\forall i \in [N] \\
& \bm{s} \geq \bm{0}, ~\bm{x} \in \mathcal{X}.
\end{array}
\end{equation}
Unfortunately, problem~\eqref{prob:individual cc reformulation} fails to be convex as its constraints involve fractions of convex functions. Below we show, however, that problem~\eqref{prob:individual cc reformulation} can be reformulated as a mixed integer conic program.  
\begin{proposition}\label{prop:individual cc}
Assume that $\bm{A}^\top\bm{x} \ne \bm{b}$ for all $\bm{x} \in \mathcal{X}$. For the safety set $\mathcal{S}(\bm{x}) = \{\bm{\xi} \in \mathbb{R}^K \mid (\bm{A}\bm{\xi} + \bm{a})^\top \bm{x} < \bm{b}^\top\bm{\xi} + b\}$, problem~\eqref{prob:cc general} is equivalent to the mixed integer conic program
\begin{equation}\label{prob:individual cc reformulation linearization}
\begin{array}{rcll}
Z^\star_{\rm ICC} =& \displaystyle \min_{\bm{q}, \bm{s}, t, \bm{x}} & \bm{c}^\top\bm{x} \\
&{\rm s.t.} & \varepsilon N t - \bm{e}^\top\bm{s} \geq \theta N \|\bm{b} - \bm{A}^\top\bm{x}\|_* \\
&& (\bm{b} - \bm{A}^\top\bm{x})^\top\bmh{\xi}_i + b - \bm{a}^\top\bm{x} + {\rm M} q_i \geq t - s_i &~\forall i \in [N] \\
&& {\rm M} (1 - q_i) \geq t - s_i &~\forall i \in [N] \\
&& \bm{q} \in \{0,1\}^N, ~\bm{s} \geq \bm{0}, ~\bm{x} \in \mathcal{X},
\end{array}
\end{equation}
where ${\rm M}$ is a suitably large (but finite) positive constant.
\end{proposition}

\noindent \emph{Proof of Proposition~\ref{prop:individual cc}.} $\;$
We already know that the chance constrained program~\eqref{prob:cc general} is equivalent to the non-convex optimization problem~\eqref{prob:individual cc reformulation}. A complicating feature of this problem is the appearance of the maximum operator in the second constraint group, which evaluates the positive part of $(\bm{b} - \bm{A}^\top\bm{x})^\top\bmh{\xi}_i + b - \bm{a}^\top\bm{x}$. To eliminate this maximum operator, for each $i\in[N]$ we introduce a binary variable $q_i \in \{0,1\}$, and we re-express the $i^{\rm th}$ member of the second constraint group via the two auxiliary constraints
\begin{equation}\label{eq:big-M}
\dfrac{(\bm{b} - \bm{A}^\top\bm{x})^\top\bmh{\xi}_i + b - \bm{a}^\top\bm{x}}{\|\bm{b} - \bm{A}^\top\bm{x}\|_*} + {\rm M} q_i \geq t - s_i 
\text{~~and~~} {\rm M} (1 - q_i) \geq t - s_i.
\end{equation}
Note that at optimality we have $q_i=1$ if $(\bm{b} - \bm{A}^\top\bm{x})^\top\bmh{\xi}_i + b - \bm{a}^\top\bm{x}$ is negative and $q_i=0$ otherwise. Intuitively, $q_i$ thus activates the less restrictive one of the two auxiliary constraints in~\eqref{eq:big-M}.
Next, we apply the variable substitutions $t\leftarrow t/\|\bm{b} - \bm{A}^\top\bm{x}\|_*$ and $\bm{s}\leftarrow \bm{s}/\|\bm{b} - \bm{A}^\top\bm{x}\|_*$, which is admissible because $\bm{A}^\top\bm{x} \ne \bm{b}$ for all $\bm{x} \in \mathcal{X}$. This change of variables yields the postulated reformulation~\eqref{prob:individual cc reformulation linearization}. 

To see that a finite value of $\rm M$ is sufficient for our reformulation to be exact, we show that the expression $((\bm{b} - \bm{A}^\top\bm{x})^\top\bmh{\xi}_i + b - \bm{a}^\top\bm{x}) / \|\bm{b} - \bm{A}^\top\bm{x}\|_*$ as well as the values of $t$ and $s_i$, $i \in [N]$, in~\eqref{eq:big-M} can all be bounded without affecting the optimal value of problem~\eqref{prob:individual cc reformulation linearization}. This is clear for the fraction as $\mathcal{X}$ is compact and the denominator is non-zero for all $\bm{x} \in \mathcal{X}$. Moreover, $t$ is nonnegative as otherwise the first constraint in~\eqref{prob:individual cc reformulation linearization} would be violated. For any fixed values of $\bm{x}$ and $t$, an optimal value of $s_i$, $i \in [N]$, is given by $s_i^\star (\bm{x}, t) = ( t -  ((\bm{b} - \bm{A}^\top\bm{x})^\top\bmh{\xi}_i + b - \bm{a}^\top\bm{x}) / \|\bm{b} - \bm{A}^\top\bm{x}\|_*)^+$. Since $\mathcal{X}$ is bounded, it thus remains to show that $t$ can be bounded from above. Indeed, for sufficiently large (but finite) $t$, the slope of $\varepsilon N t - \bm{e}^\top \bm{s}^\star (\bm{x}, t)$ on the left-hand side of the first constraint in~\eqref{prob:individual cc reformulation linearization} is $- (1 - \varepsilon) N$. Since $\varepsilon < 1$, we thus conclude that this constraint is violated for large values of $t$.
\hfill \Halmos
\endproof

\begin{remark}
The condition that $\bm{A}^\top\bm{x} \ne \bm{b}$ for all $\bm{x} \in \mathcal{X}$ does not restrict the generality of our formulation. Indeed, if an optimal solution $(\bm{q}^\star, \bm{s}^\star, t^\star, \bm{x}^\star)$ to problem~\eqref{prob:individual cc reformulation linearization} satisfies $\bm{A}^\top\bm{x}^\star \ne \bm{b}$, then $\bm{x}^\star$ solves problem~\eqref{prob:cc general} since our argument in the proof of Proposition~\ref{prop:individual cc} applies to $\bm{x}^\star$ even if $\bm{A}^\top\bm{x} = \bm{b}$ for some $\bm{x} \in \mathcal{X}$. Assume now that an optimal solution $(\bm{q}^\star, \bm{s}^\star, t^\star, \bm{x}^\star)$ to problem~\eqref{prob:individual cc reformulation linearization} satisfies $\bm{A}^\top\bm{x}^\star = \bm{b}$. In that case, the ambiguous chance constraint in problem~\eqref{prob:cc general} requires that $\bm{a}^\top \bm{x}^\star < b$. If that is the case for $\bm{x}^\star$, it is optimal in problem~\eqref{prob:cc general}. If, finally, an optimal solution $(\bm{q}^\star, \bm{s}^\star, t^\star, \bm{x}^\star)$ to problem~\eqref{prob:individual cc reformulation linearization} satisfies $\bm{A}^\top\bm{x}^\star = \bm{b}$ and $\bm{a}^\top \bm{x}^\star \geq b$, then one would ideally like to solve a variant of problem~\eqref{prob:individual cc reformulation linearization} that includes the additional constraint
\begin{equation}\label{rem2_constraint}
\bm{A}^\top \bm{x} \neq \bm{b} \quad \text{or} \quad \bm{a}^\top \bm{x} < b.
\end{equation}
This variant of problem~\eqref{prob:individual cc reformulation linearization} can be solved by solving $2 K + 1$ versions of problem~\eqref{prob:individual cc reformulation linearization}, where each version includes exactly one of the constraints $[\bm{A}^\top\bm{x}]_k > [\bm{b}]_k$, $[\bm{A}^\top\bm{x}]_k < [\bm{b}]_k$, $k \in [K]$, or $\bm{a}^\top \bm{x} < b$. One readily verifies that the solution that attains the least objective value amongst these $2K + 1$ versions of problem~\eqref{prob:individual cc reformulation linearization} is an optimal solution to problem~\eqref{prob:individual cc reformulation linearization} with the added constraint~\eqref{rem2_constraint}.
\end{remark}

\begin{remark}
The mixed-integer conic program~\eqref{prob:individual cc reformulation linearization} simplifies to a mixed-integer linear program whenever $\|\cdot\|$ represents the $1$-norm or the $\infty$-norm, and it can be reformulated as a mixed-integer second-order cone program whenever $\|\cdot\|$ represents a $p$-norm for some $p \in \mathbb{Q}$, $p>1$, see Section~2.3.1 in \cite{Ben-tal_Nemirovski_book}.
\end{remark}

\begin{remark}
The deterministic reformulation~\eqref{prob:individual cc reformulation linearization} is remarkably parsimonious. For an $L$-dimensional feasible region $\mathcal{X} \subseteq \mathbb{R}^L$ and an empirical distribution $\hat{\mathbb{P}}$ with $N$ data points, our reformulation~\eqref{prob:individual cc reformulation linearization} has $N$ binary variables, $L + N + 1$ continuous decisions as well as $2N + 1$ constraints (excluding those that describe $\mathcal{X}$). In comparison, a classical chance constrained formulation, which is tantamount to setting the Wasserstein radius to $\theta = 0$ in problem~\eqref{prob:cc general}, has $N$ binary variables, $L$ continuous decisions as well as $N + 1$ constraints. Thus, adding distributional robustness only requires an additional $N + 1$ continuous decisions as well as $N$ further constraints.
\end{remark}

\begin{remark}\label{rem:piecewise_linear}
The deterministic reformulation~\eqref{prob:individual cc reformulation linearization} requires the specification of a sufficiently large constant $\mathrm{M}$, which can typically be determined by an investigation of the structure of problem~\eqref{prob:individual cc reformulation linearization}. Alternatively, many commercial solver packages allow to directly specify the following reformulation of problem~\eqref{prob:individual cc reformulation linearization} via the use of piecewise linear constraints:
\begin{equation*}
\begin{array}{rcll}
Z^\star_{\rm ICC} =& \displaystyle \min_{\bm{q}, \bm{s}, t, \bm{x}} & \bm{c}^\top\bm{x} \\
&{\rm s.t.} & \varepsilon N t - \bm{e}^\top\bm{s} \geq \theta N \|\bm{b} - \bm{A}^\top\bm{x}\|_* \\
&& ( (\bm{b} - \bm{A}^\top\bm{x})^\top\bmh{\xi}_i + b - \bm{a}^\top\bm{x})^+ \geq t - s_i &~\forall i \in [N] \\
&& \bm{s} \geq \bm{0}, ~\bm{x} \in \mathcal{X}
\end{array}
\end{equation*}
This formulation has the advantage that it does not require the specification of the constant $\mathrm{M}$.
\end{remark}

\subsection{Reformulation of Joint Chance Constraints with Right-Hand Side Uncertainty}\label{sec:ref_joint_cc}

Assume next that problem~\eqref{prob:cc general} accommodates a joint chance constraint defined through the safety set $\mathcal{S}(\bm{x}) = \{\bm{\xi} \in \mathbb{R}^K \mid \bm{a}^\top_m \bm{x} < \bm{b}^\top_m\bm{\xi} + b_m ~\forall m \in [M]\}$, in which the uncertainty affects only the right-hand sides of the safety conditions. Without loss of generality, we may assume that $\bm{b}_m \ne \bm{0}$ for all $m \in [M]$. Indeed, if $\bm{b}_m = \bm{0}$, then the $m^{\rm th}$ safety condition in the chance constraint becomes independent of the uncertainty and can thus be absorbed in $\mathcal{X}$. Joint chance constrained programs with right-hand side uncertainty have been proposed, among others, for problems in transportation \citep{luedtke2010}, lot-sizing \citep{beraldi2002branch, kuccukyavuz2012mixing}, unit commitment \citep{yanagisawa_2013} and project management \citep{wiesemann2012multi}.

Observe that the complement of the safety set is now representable as $\bar{\mathcal{S}}(\bm{x}) = \bigcup_{m \in [M]} \mathcal{H}_m(\bm{x})$, where $\mathcal{H}_m(\bm{x}) = \{\bm{\xi} \in \mathbb{R}^K \mid \bm{a}^\top_m \bm{x} \geq \bm{b}^\top_m\bm{\xi} + b_m\}$ is a closed halfspace for every $m \in [M]$. By Lemma~\ref{lem:distance to the union of closed half-spaces} in the appendix we have
\begin{equation}\label{eq:distance-joint}
\mathbf{dist}(\bmh{\xi}_i, \bar{\mathcal{S}}(\bm{x})) = \min_{m \in [M]} \bigg\{ \dfrac{(\bm{b}^\top_m\bmh{\xi}_i + b_m - \bm{a}^\top_m\bm{x})^+}{\|\bm{b}_m\|_*} \bigg\} = \bigg(\min_{m \in [M]} \bigg\{\dfrac{\bm{b}^\top_m\bmh{\xi}_i + b_m - \bm{a}^\top_m\bm{x}}{\|\bm{b}_m\|_*} \bigg\}\bigg)^+.
\end{equation}
With this closed-form expression for the distance to the unsafe set, we can reformulate problem~\eqref{prob:cc general} as a mixed integer conic program.

\begin{proposition}\label{prop:joint cc}
For the safety set $\mathcal{S}(\bm{x}) = \{\bm{\xi} \in \mathbb{R}^K \mid \bm{a}^\top_m \bm{x} < \bm{b}^\top_m\bm{\xi} + b_m ~\forall m \in [M]\}$, where $\bm{b}_m \ne \bm{0}$ for all $m \in [M]$, the chance constrained program~\eqref{prob:cc general} is equivalent to the mixed integer conic program
\begin{equation}
\label{prob:joint cc reformulation M linearization}
\begin{array}{rcll}
Z^\star_{\rm JCC} =& \displaystyle \min_{\bm{q}, \bm{s}, t, \bm{x}} & \bm{c}^\top\bm{x} \\
&{\rm s.t.} & \varepsilon N t - \bm{e}^\top\bm{s} \geq \theta N \\
&& \dfrac{\bm{b}^\top_m\bmh{\xi}_i + b_m - \bm{a}^\top_m\bm{x}}{\|\bm{b}_m\|_*} + {\rm M} q_i \geq t - s_i &~\forall m \in [M], ~i \in [N] \\
&& {\rm M} (1 - q_i) \geq t - s_i &~\forall i \in [N] \\
&& \bm{q} \in \{0,1\}^N,~ \bm{s} \geq \bm{0}, ~\bm{x} \in \mathcal{X},
\end{array}
\end{equation}
where ${\rm M}$ is a suitably large (but finite) positive constant.
\end{proposition}

\noindent \emph{Proof of Proposition~\ref{prop:joint cc}.} $\;$
By Theorem~\ref{thm:cc-deterministic}, the chance constrained program~\eqref{prob:cc general} is equivalent to~\eqref{prob:cc reformulation}. Using~\eqref{eq:distance-joint}, the $i^{\rm th}$~member of the second constraint group in~\eqref{prob:cc reformulation} can be reformulated as
$$
\bigg(\displaystyle \min_{m \in [M]} \bigg\{\dfrac{\bm{b}^\top_m\bmh{\xi}_i + b_m - \bm{a}^\top_m\bm{x}}{\|\bm{b}_m\|_*} \bigg\}\bigg)^+ \geq t - s_i.
$$
To eliminate the maximum operator, we introduce a binary variable $q_i \in \{0,1\}$ to re-express the above constraint as
$$
\left\{
\begin{array}{ll}
\dfrac{\bm{b}^\top_m\bmh{\xi}_i + b_m - \bm{a}^\top_m\bm{x}}{\|\bm{b}_m\|_*} + {\rm M} q_i \geq t - s_i &~\forall m \in [M]\\
{\rm M} (1 - q_i) \geq t - s_i 
\end{array}
\right.
$$
A similar argument as in the proof of Proposition~\ref{prop:individual cc} shows that a finite value of $\rm M$ is sufficient for our reformulation to be exact.
\hfill \Halmos
\endproof

Similar to Remark~\ref{rem:piecewise_linear} in the previous section, many commercial solvers allow to directly specify a reformulation of  problem~\eqref{prob:joint cc reformulation M linearization} that replaces the constant $\mathrm{M}$ with piecewise linear constraints.

\begin{remark}
The deterministic reformulation~\eqref{prob:joint cc reformulation M linearization} has $N$ binary variables, $L + N + 1$ continuous decisions as well as $(M + 1) N + 1$ constraints (excluding those that describe $\mathcal{X}$). In comparison, the corresponding classical chance constrained formulation has $N$ binary variables, $L$ continuous decisions as well as $MN + 1$ constraints. Thus, adding distributional robustness requires an additional $N + 1$ continuous decisions as well as $N$ further (linear) constraints.
\end{remark}

\section{Numerical Experiments}

We compare our exact reformulation of the ambiguous chance constrained program~\eqref{prob:cc general} with the bicriteria approximation scheme of \cite{xie2020bicriteria} on a portfolio optimization problem in Section~\ref{sec:portfolio} as well as with a classical (non-ambiguous) chance constrained formulation and a Kernel density estimator based version of the ambiguous chance constrained program over a $\phi$-divergence ambiguity set on a transportation problem in Section~\ref{sec:transportation}. Our goal is to investigate the computational scalability of our reformulation as well as its out-of-sample performance in a data-driven setting. All results were produced on an Intel Xeon 2.66GHz processor with 8GB memory in single-core mode using CPLEX 12.8. Following Remark~\ref{rem:piecewise_linear}, we avoid the specification of the constant $\mathrm{M}$ in our ambiguous chance constrained program through the use of piecewise linear constraints.

\subsection{Portfolio Optimization}\label{sec:portfolio}

We consider a portfolio optimization problem studied by \cite{xie2020bicriteria}. The problem asks for the minimum-cost portfolio investment $\bm{x}$ into $K$ assets with random returns $\tilde{\xi}_1, \ldots, \tilde{\xi}_K$ that exceeds a pre-specified target return $w$ with high probability $1 - \varepsilon$. The problem can be cast as the following instance of the ambiguous chance constrained program~\eqref{prob:cc general}:
\begin{equation}\label{eq:portfolio}
\begin{array}{cll}
\displaystyle \min_{\bm{x}} & \bm{c}^\top\bm{x} \\
{\rm s.t.} &\displaystyle \mathbb{P}[\bmt{\xi}^\top\bm{x} > w] \geq 1-\varepsilon &~\forall \mathbb{P} \in \mathcal{F}(\theta)\\
&\bm{x} \geq \bm{0}.
\end{array}
\end{equation}

We compare our exact reformulation of problem~\eqref{eq:portfolio} with the $(\sigma, \gamma)$-bicriteria approximation scheme of \cite{xie2020bicriteria}, which produces solutions that satisfy the ambiguous chance constraint in~\eqref{eq:portfolio} with probability $1 - \sigma \varepsilon$, $\sigma > 1$, and whose costs are guaranteed to exceed the optimal costs in~\eqref{eq:portfolio} by a factor of at most $\gamma = \sigma / (\sigma - 1)$. Since the bicriteria approximation scheme can readily utilize support information for the random vector $\bmt{\xi}$, we replace the ambiguity set $\mathcal{F}(\theta)$ with $\bar{\mathcal{F}}(\theta) = \mathcal{F}(\theta) \cap \{\mathbb{P} \mid \mathbb{P}[\bmt{\xi} \in \mathbb{R}^K_+] = 1\}$ in their approach. Contrary to the experiments conducted by \cite{xie2020bicriteria}, we set $\sigma = 1$. This is to the disadvantage of their approach, as it does not provide any approximation guarantees in that case, but it allows us to compare the resulting portfolios as they provide the same return guarantees. For the performance of the bicriteria approximation scheme with $\sigma > 1$, we refer to Section~6.2 of \cite{xie2020bicriteria}.

In our numerical experiments, we consider a similar setting as \cite{xie2020bicriteria}. We set $K = 50$, $w = 1$ and choose the cost coefficients $c_1, \ldots, c_{50}$ uniformly at random from $\{ 1, \ldots, 100 \}$. Each asset return $\tilde{\xi}_i$ is governed by a uniform distribution on $[0.8, 1.5]$, and we assume that $N = 100$ training samples $\bmh{\xi}_1, \ldots, \bmh{\xi}_{100}$ are available. We use the $2$-norm Wasserstein ambiguity set, which implies that our exact reformulation of problem~\eqref{eq:portfolio} is a mixed-integer second-order cone program, and set the Wasserstein radius to $\theta \in  \{0.05, 0.1, 0.2\}$. The risk threshold is set to $\varepsilon \in \{0.05, 0.1\}$.

\begin{table}[tb]
\begin{center}
\begin{tabular}{c>{\centering\arraybackslash}p{1.5cm}>{\centering\arraybackslash}>{\centering\arraybackslash}p{1.5cm}>{\centering\arraybackslash}>{\centering\arraybackslash}p{1.5cm}>{\centering\arraybackslash}p{1.5cm}>{\centering\arraybackslash}p{1.5cm}>{\centering\arraybackslash}p{1.5cm}}
\hline
\hline  
\multirow{2}[0]{*}{$(\varepsilon, \theta)$}&\multicolumn{3}{c}{Ratio of objective values} &\multicolumn{3}{c}{Ratio of runtimes} \\
\cline{2-7}
&5\%  &50\%  &95\% &5\%  &50\%  &95\% \\ 
\hline
\multicolumn{1}{c}{$(0.05, 0.05)$} &1.6 &2.4 &3.2 &5.2 &8.3 &10.8  \\ 
\multicolumn{1}{c}{$(0.05, 0.10)$} &1.9 &2.9 &5.0 &4.9 &7.7 &10.6 \\ 
\multicolumn{1}{c}{$(0.05, 0.20)$} &2.3 &2.8 &3.5 &3.8 &4.9 &7.2  \\ 
\multicolumn{1}{c}{$(0.10, 0.05)$} &1.0 &1.1 &1.3 &7.3 &10.9 &13.0  \\ 
\multicolumn{1}{c}{$(0.10, 0.10)$} &1.5 &2.3 &3.1 &7.1 &9.7 &13.3  \\
\multicolumn{1}{c}{$(0.10, 0.20)$} &2.1 &2.7 &3.9 &4.2 &6.2 &10.1  \\ 
\hline
\hline \\
\end{tabular}
\end{center}
\caption{\textnormal{Objective and runtime ratios of the bicriteria approximation scheme for different values of $\varepsilon$ and $\theta$. For each parameter setting, we report the $5\%$, $50\%$ and $95\%$ quantiles over 50 randomly generated instances.}}
\label{table:bicteria}
\end{table}

\begin{figure}[tb]
\begin{subfigure}{.5\textwidth}
\begin{center}
\includegraphics[width=1\linewidth]{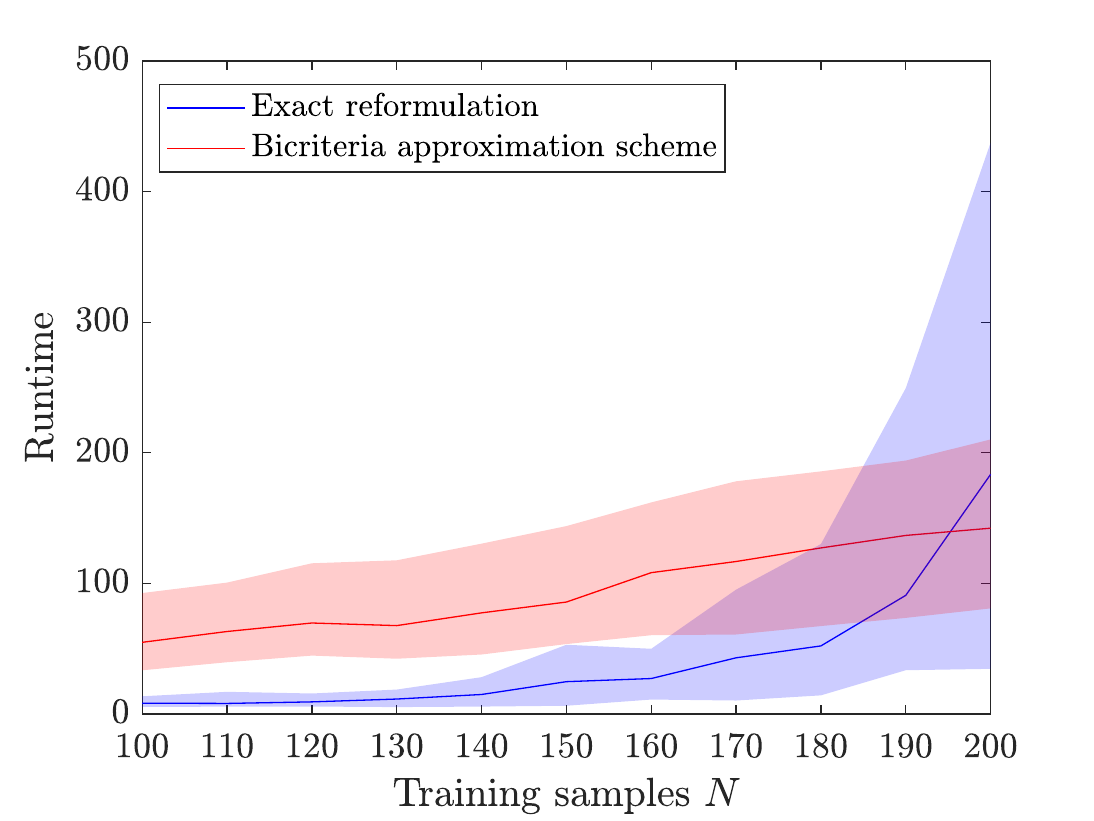}
\end{center}
\end{subfigure}
\begin{subfigure}{.5\textwidth}
\begin{center}
\includegraphics[width=1\linewidth]{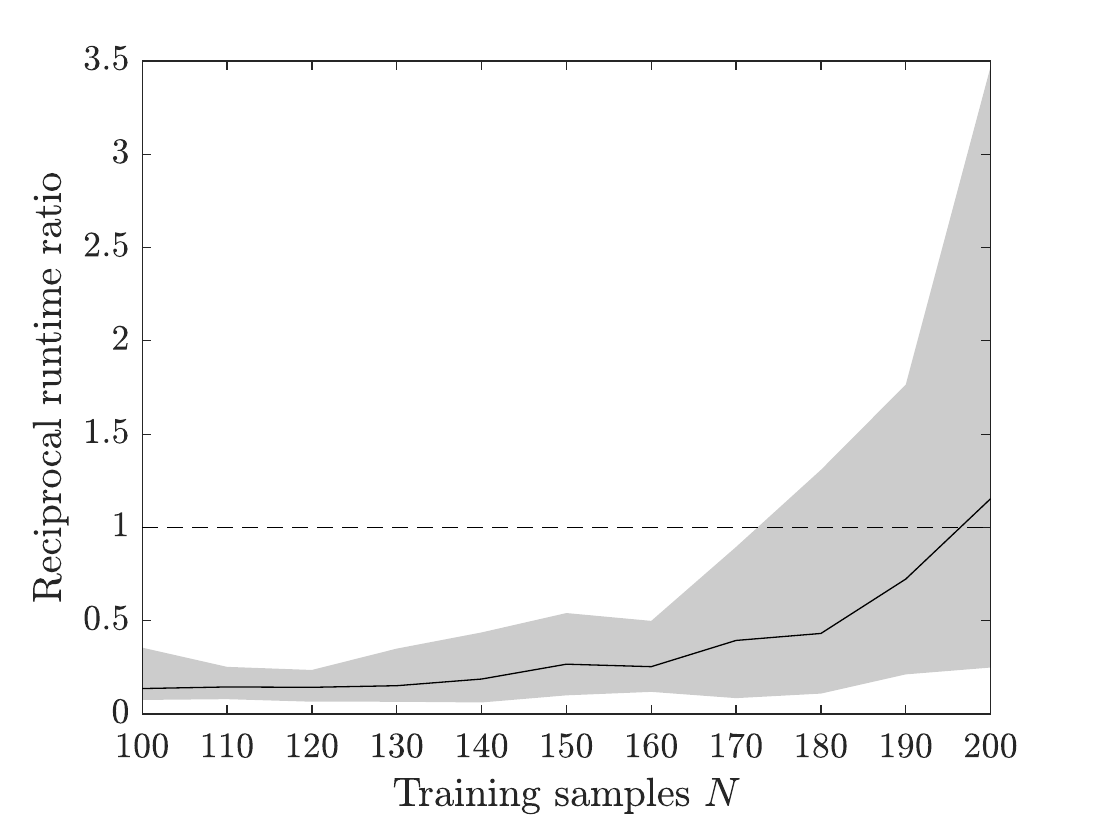}
\end{center}
\end{subfigure}
\vspace{0.2cm}
\caption{{\textnormal{Runtimes (left) and reciprocal runtime ratios (right) of our exact reformulation and the bicriteria approximation scheme for $(\varepsilon, \theta) = (0.10,0.05)$ and different sample sizes $N$. The shaded regions cover the $5\%$ to $95\%$ quantiles of $50$ randomly generated instances, whereas the solid lines describe the median statistics.}} \label{fig:runtime}}
\end{figure}

Table~\ref{table:bicteria} compares the objective values and runtimes of our exact reformulation and the bicriteria approximation scheme for various combinations of the risk threshold $\varepsilon$ and Wasserstein radius $\theta$. The table shows that despite incorporating additional support information, the bicriteria approximation scheme determines solutions whose costs significantly exceed those of the solutions found by our exact reformulation. Perhaps more surprisingly, the bicriteria approximation scheme is also computationally more expensive. As Figure~\ref{fig:runtime} shows, however, this is an artifact of the small sample size $N$ employed in the experiments of \cite{xie2020bicriteria}, and the bicriteria approximation scheme is faster than our exact reformulation for larger samples sizes.

\subsection{Transportation}\label{sec:transportation}
We consider a probabilistic transportation problem studied by \cite{luedtke2010} and \cite{yanagisawa_2013}. The problem asks for the cost-optimal distribution of a single good from a set of factories $f \in [F]$ to a set of distribution centers $d \in [D]$. Each factory $f \in [F]$ has an individual production capacity $m_f$, and each distribution center $d \in [D]$ faces a random aggregate customer demand $\tilde{\xi}_d$. The cost of shipping one unit of the good from factory $f$ to distribution center $d$ is denoted by $c_{fd}$. We aim to find a transportation plan that minimizes the shipping costs, respects the production capacity of each factory and satisfies the demand at each distribution center with high probability. The problem can be cast as the following instance of problem~\eqref{prob:cc general}:
\begin{equation}\label{eq:amb_transp_prob}
\begin{array}{cl@{\quad}l}
\displaystyle \min_{\bm{x}} & \bm{c}^\top\bm{x} \\
{\rm s.t.} & \displaystyle \mathbb{P} \Bigg[ \sum_{f \in [F]} x_{fd} \geq \tilde{\xi}_d \quad \forall d \in [D] \Bigg] \geq 1 - \varepsilon &~ \forall \mathbb{P} \in \mathcal{F}(\theta) \\[5mm]
& \displaystyle \sum_{d \in [D]} x_{fd} \leq m_f &~ \forall f \in [F] \\
& \displaystyle \bm{x} \geq \bm{0}.
\end{array}
\end{equation}
Here, $x_{fd}$ denotes the quantity shipped from factory $f \in [F]$ to distribution center $d \in [D]$. Problem~\eqref{eq:amb_transp_prob} is an ambiguous joint chance constrained program with right-hand side uncertainty. Since each safety condition in~\eqref{eq:amb_transp_prob} contains a single random variable with coefficient $1$ on the right-hand side, our exact reformulation reduces to the same mixed-integer linear program for any norm $\lVert \cdot \rVert$.

In our first experiment, we investigate the scalability of the exact reformulation of problem~\eqref{eq:amb_transp_prob} that is offered by Proposition~\ref{prop:joint cc}. To this end, we generate random test instances with $5$ factories and $10, 20, \ldots, 50$ distribution centers that are located uniformly at random on the Euclidean plane $[0, 10]^2$. We identify the transportation costs $c_{fd}$ with the Euclidean distances between the factories and distribution centers. The demand vector $\bmt{\xi}$ is described by $50$, $100$ or $150$ samples from a uniform distribution that is supported on $[0.8 \bm{\mu}, 1.2 \bm{\mu}]$, where the expected demand $\mu_d$ at distribution center $d \in [D]$ is picked uniformly at random from the interval $[0, 10]$. The capacity of each factory is chosen uniformly at random, and the capacities are subsequently scaled so that the factories can jointly produce up to $150\%$ of the maximum cumulative demand. For each instance, we choose $10$ ascending Wasserstein radii $\theta_1 < \ldots < \theta_{10}$ uniformly so that $\theta_1 = 0.001$ and $\theta_{10}$ is the smallest radius for which the corresponding instance of problem~\eqref{eq:amb_transp_prob} becomes infeasible. We fix $\varepsilon = 0.1$.

\begin{figure}[tb]
\begin{subfigure}{.33\textwidth}
\begin{center}
\includegraphics[width=0.9\linewidth]{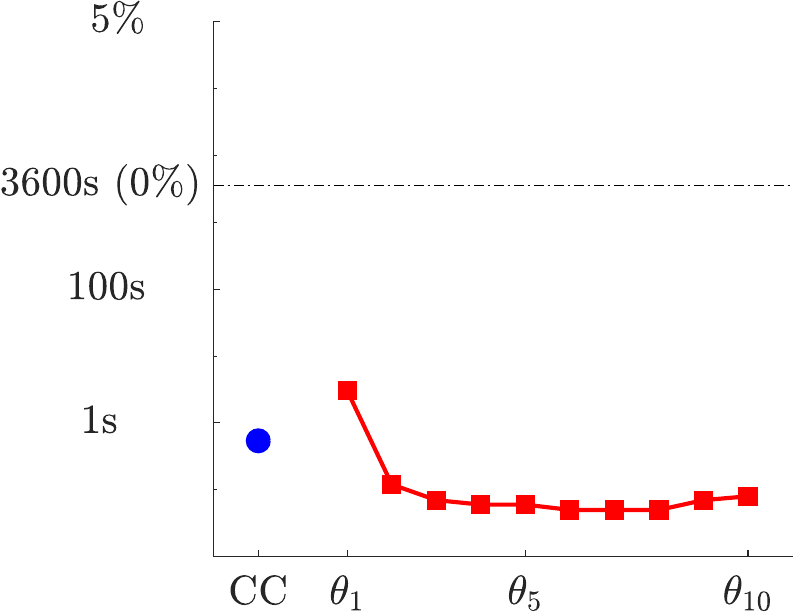}
\end{center}
\end{subfigure}%
\begin{subfigure}{.33\textwidth}
\begin{center}
\includegraphics[width=0.9\linewidth]{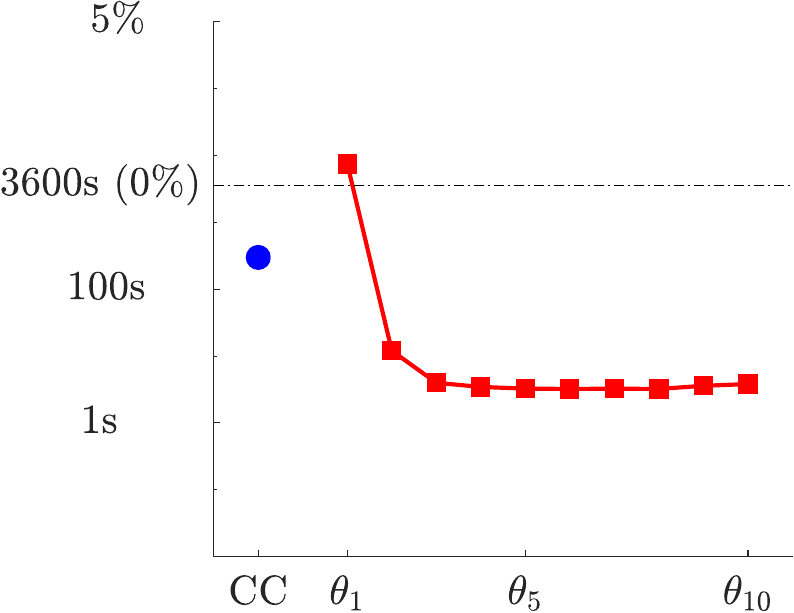}
\end{center}
\end{subfigure}	
\begin{subfigure}{.33\textwidth}
\begin{center}
\includegraphics[width=0.9\linewidth]{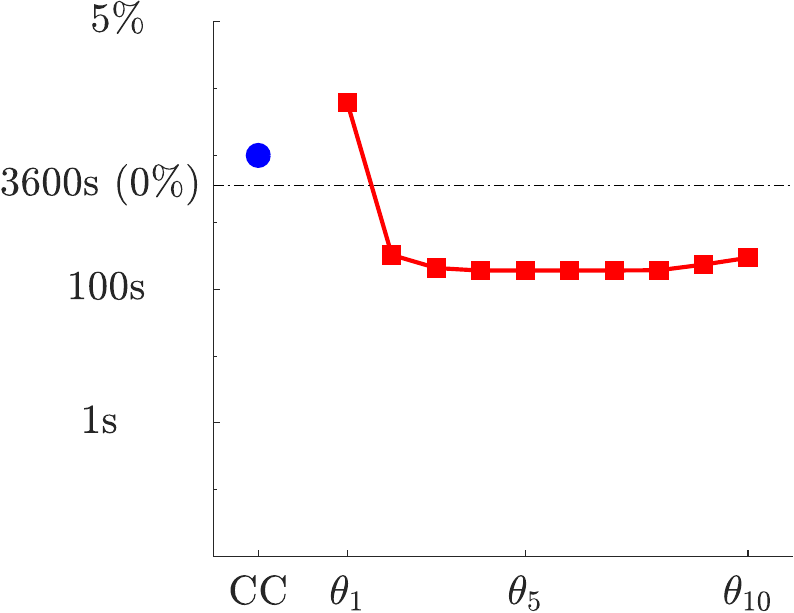}
\end{center}
\end{subfigure}
\vspace{0.2cm}
\caption{{\textnormal{Median solution times (below dashed lines) and optimality gaps (above dashed lines) for $D = 10$ and $N = 50$ (left), $D = 30$ and $N = 100$ (middle) and $D = 50$ and $N = 150$ (right).}} \label{fig:case}}
\end{figure}

\begin{figure}[tb]
\begin{subfigure}{.33\textwidth}
\begin{center}
\includegraphics[width=0.95\linewidth]{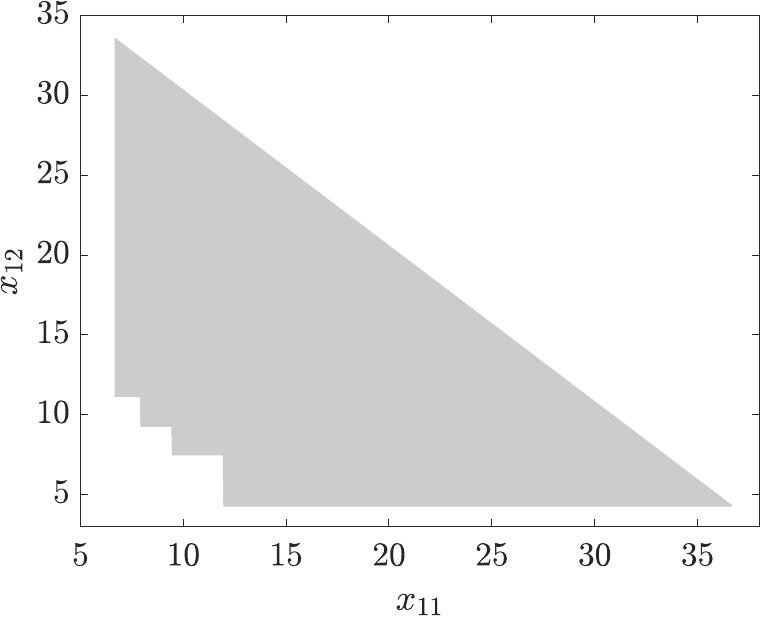}
\end{center}
\end{subfigure}%
\begin{subfigure}{.33\textwidth}
\begin{center}
\includegraphics[width=0.95\linewidth]{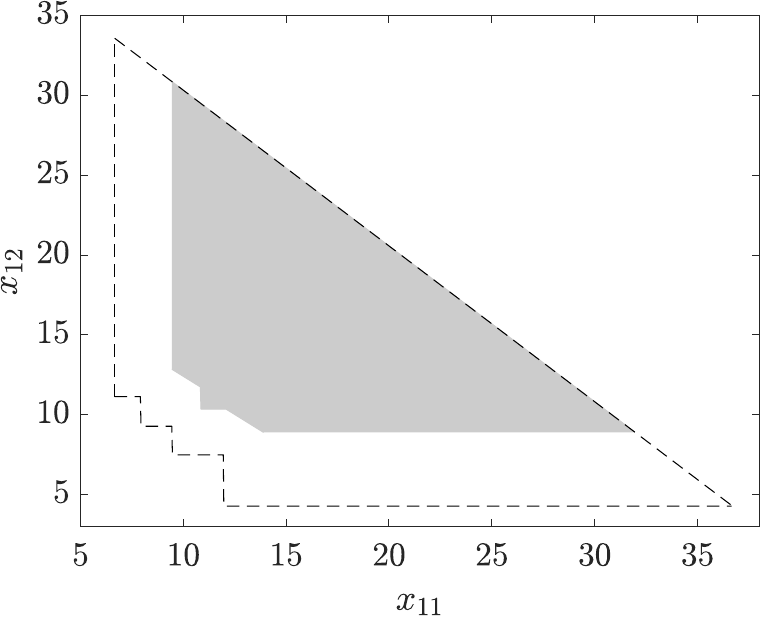}
\end{center}
\end{subfigure}	
\begin{subfigure}{.33\textwidth}
\begin{center}
\includegraphics[width=0.95\linewidth]{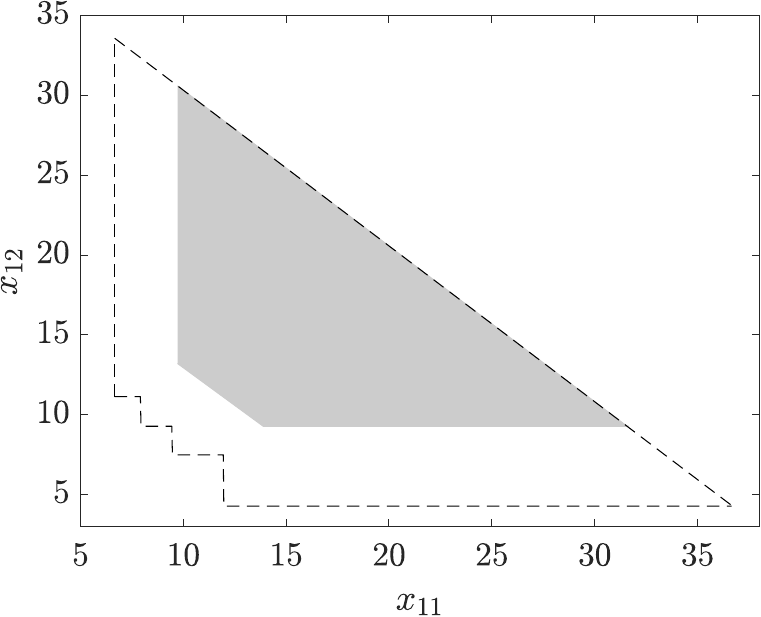}
\end{center}
\end{subfigure}
\vspace{0.2cm}
\caption{{\textnormal{For a transportation problem with $F = 1$ factory, $D = 2$ distribution centers and $N = 10$ training samples, the graphs visualize the feasible regions of the classical chance constrained formulation~\eqref{eq:classical_transp_prob} (left) and the ambiguous chance constrained problem~\eqref{eq:amb_transp_prob} for a small (middle) and a large (right) value of $\theta$.}} \label{fig:radius}}
\end{figure}

Tables~\ref{table:scalability_50}--\ref{table:scalability_150} and Figure~\ref{fig:case} compare the runtimes of our ambiguous chance constrained program with those of the classical chance constrained formulation of problem~\eqref{eq:amb_transp_prob},
\begin{equation}\label{eq:classical_transp_prob}
\begin{array}{cl@{\quad}l}
\displaystyle \min_{\bm{x}, \bm{y}} & \bm{c}^\top\bm{x} \\
{\rm s.t.} & \displaystyle \sum_{f \in [F]} x_{fd} + \mathrm{M} y_i \geq \hat{\xi}_{id} & \forall d \in [D], ~i \in [N] \\
& \displaystyle \bm{e}^\top \bm{y} \leq \lfloor \varepsilon N \rfloor \\
& \displaystyle \sum_{d \in [D]} x_{fd} \leq m_f & \displaystyle \forall f \in [F] \\
& \displaystyle \bm{x} \geq \bm{0}, ~ \bm{y} \in \{ 0, 1 \}^N,
\end{array}
\end{equation}
where $\mathrm{M}$ is a sufficiently large positive constant. The results show that for the smallest Wasserstein radius $\theta_1 = 0.001$, the ambiguous chance constrained program~\eqref{eq:amb_transp_prob} is---as expected---more difficult to solve than the corresponding classical chance constrained program~\eqref{eq:classical_transp_prob}. Interestingly, the ambiguous chance constrained program becomes considerably \emph{easier} to solve than the classical chance constrained program for the larger Wasserstein radii $\theta_2, \ldots, \theta_{10}$. This surprising result is explained in Figure~\ref{fig:radius}, which shows that the feasible region of the ambiguous chance constrained program tends to convexify as the Wasserstein radius $\theta$ increases. In fact, one can show that the set of vectors $\bm{q} \in \{ 0, 1 \}^N$ that are feasible in the deterministic reformulation of problem~\eqref{eq:amb_transp_prob} shrinks monotonically with $\theta$. Since it is the presence of these binary vectors that causes the non-convexity of problem~\eqref{eq:amb_transp_prob}, one can expect the problem to become better behaved as $\theta$ increases.

\begin{table}[tb]
\begin{center}
\begin{tabular}{c>{\centering\arraybackslash}p{1.1cm}>{\centering\arraybackslash}p{0.9cm}>{\centering\arraybackslash}p{0.9cm}>{\centering\arraybackslash}p{0.9cm}>{\centering\arraybackslash}p{0.9cm}>{\centering\arraybackslash}p{0.9cm}>{\centering\arraybackslash}p{0.9cm}>{\centering\arraybackslash}p{0.9cm}>{\centering\arraybackslash}p{0.9cm}>{\centering\arraybackslash}p{0.9cm}>{\centering\arraybackslash}p{0.9cm}>{\centering\arraybackslash}p{0.9cm}}\hline
\hline
{\tabincell{c}{$\#$ of distribution \\ centers}} & CC  &$\theta_1$  &$\theta_2$ &$\theta_3$ &$\theta_4$  &$\theta_5$
&$\theta_6$  &$\theta_7$
&$\theta_8$  &$\theta_9$ &$\theta_{10}$ \\ 
\hline
10 &0.5 &3.0 &0.1 &$<0.1$ &$<0.1$ &$<0.1$ &$<0.1$ &$<0.1$ &$<0.1$ &$<0.1$ &$<0.1$ \\ 
20 &4.0 &9.7 &0.2 &0.1 &0.1 &0.1 &$<0.1$ &$<0.1$ &$<0.1$ &0.1 &0.1 \\ 
30 &7.3 &13.1 &0.3 &0.2 &0.1 &0.1 &0.1 &0.1 &0.1 &0.1 &0.2 \\ 
40 &11.2 &19.3 &0.4 &0.2 &0.2 &0.2 &0.2 &0.2 &0.2 &0.2 &0.3 \\ 
50 &15.8 &166.5 &0.3 &0.2 &0.2 &0.2 &0.2 &0.2 &0.2 &0.3 &0.3 \\
\hline
\hline \\
\end{tabular}
\end{center}
\caption{\textnormal{Solution times in seconds for $N = 50$ training samples. `CC' and `$\theta_i$' refer to problem~\eqref{eq:classical_transp_prob} and problem~\eqref{eq:amb_transp_prob} with different Wasserstein radii, respectively. We present median results over 100 random instances. Where the median solution time exceeds 3,600s, we report the median optimality gap in brackets.}}
\label{table:scalability_50}
\end{table}

\begin{table}[tb]
\begin{center}
\begin{tabular}{c>{\centering\arraybackslash}p{1.1cm}>{\centering\arraybackslash}p{0.9cm}>{\centering\arraybackslash}p{0.9cm}>{\centering\arraybackslash}p{0.9cm}>{\centering\arraybackslash}p{0.9cm}>{\centering\arraybackslash}p{0.9cm}>{\centering\arraybackslash}p{0.9cm}>{\centering\arraybackslash}p{0.9cm}>{\centering\arraybackslash}p{0.9cm}>{\centering\arraybackslash}p{0.9cm}>{\centering\arraybackslash}p{0.9cm}>{\centering\arraybackslash}p{0.9cm}}\hline
\hline
{\tabincell{c}{$\#$ of distribution \\ centers}} & CC  &$\theta_1$  &$\theta_2$ &$\theta_3$ &$\theta_4$  &$\theta_5$
&$\theta_6$  &$\theta_7$
&$\theta_8$  &$\theta_9$ &$\theta_{10}$ \\ 
\hline
10 &16.3 &166.4 &4.7 &2.0 &1.5 &1.4 &1.4 &1.4 &1.4 &1.5 &1.8 \\ 
20 &93.6 &1910.8 &8.1 &2.9 &2.5 &2.5 &2.4 &2.4 &2.4 &2.7 &2.8  \\ 
30 &298.3 &$[0.2\%]$ &12.0 &4.0 &3.5 &3.3 &3.2 &3.3 &3.2 &3.6 &3.8 \\ 
40 &664.2 &$[0.8\%]$ &16.0 &5.1 &4.7 &4.5 &4.5 &4.5 &4.4 &4.8 &5.1 \\ 
50 &1,293.2 &$[0.8\%]$ &20.3 &6.5 &5.6 &5.5 &5.4 &5.4 &5.4 &5.7 &6.2 \\
\hline
\hline \\
\end{tabular}
\end{center}
\caption{\textnormal{Solution times for $N = 100$ training samples. The table has the same interpretation as Table~\ref{table:scalability_50}.}}
\label{table:scalability_100}
\end{table}

\begin{table}[tb]
\begin{center}
\begin{tabular}{c>{\centering\arraybackslash}p{1.1cm}>{\centering\arraybackslash}p{0.9cm}>{\centering\arraybackslash}p{0.9cm}>{\centering\arraybackslash}p{0.9cm}>{\centering\arraybackslash}p{0.9cm}>{\centering\arraybackslash}p{0.9cm}>{\centering\arraybackslash}p{0.9cm}>{\centering\arraybackslash}p{0.9cm}>{\centering\arraybackslash}p{0.9cm}>{\centering\arraybackslash}p{0.9cm}>{\centering\arraybackslash}p{0.9cm}>{\centering\arraybackslash}p{0.9cm}}\hline
\hline
{\tabincell{c}{$\#$ of distribution \\ centers}} & CC  &$\theta_1$  &$\theta_2$ &$\theta_3$ &$\theta_4$  &$\theta_5$
&$\theta_6$  &$\theta_7$
&$\theta_8$  &$\theta_9$ &$\theta_{10}$ \\ 
\hline
10 &94.6 &$[0.7\%]$ &85.6 &48.5 &44.8 &44.0 &42.5 &43.3 &43.0 &52.0 &77.0 \\ 
20 &874.2 &$[1.9\%]$ &143.9 &90.5 &76.3 &75.6 &72.8 &72.5 &73.2 &85.7 &112.4 \\ 
30 &$[0.1\%]$ &$[3.2\%]$ &213.8 &126.4 &113.0 &109.5 &108.9 &108.8 &110.3 &125.4 &165.1 \\ 
40 &$[0.3\%]$ &$[3.7\%]$ &286.8 &168.2 &154.2 &149.1 &149.3 &151.7 &152.1 &182.8 &231.5 \\ 
50 &$[0.4\%]$ &$[3.0\%]$ &324.6 &207.0 &189.3 &190.9 &190.0 &190.4 &191.8 &233.0 &294.4\\
\hline
\hline \\
\end{tabular}
\end{center}
\caption{\textnormal{Solution times for $N = 150$ training samples. The table has the same interpretation as Table~\ref{table:scalability_50}.}}
\label{table:scalability_150}
\end{table}

\begin{figure}[tb]
\begin{subfigure}{.5\textwidth}
\begin{center}
\includegraphics[width=1\linewidth]{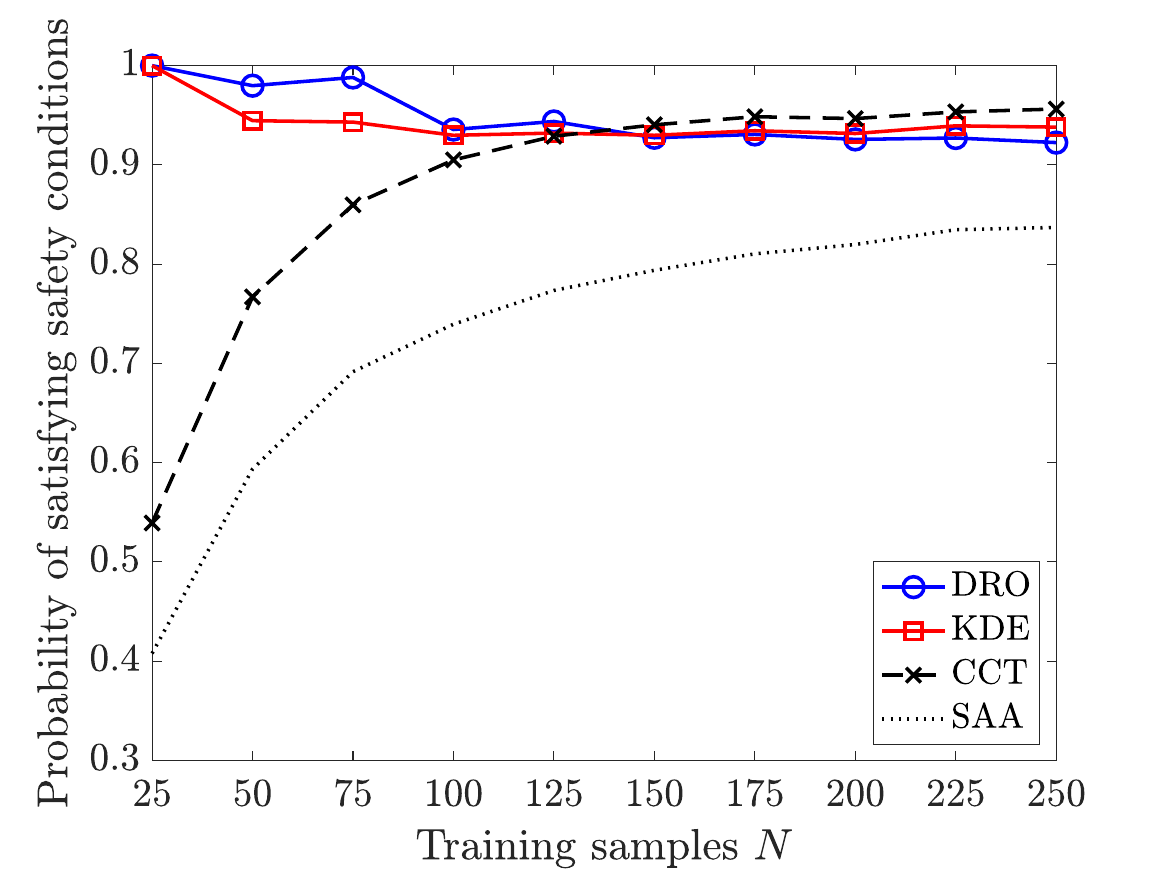}
\end{center}
\end{subfigure}%
\begin{subfigure}{.5\textwidth}
\begin{center}
\includegraphics[width=1\linewidth]{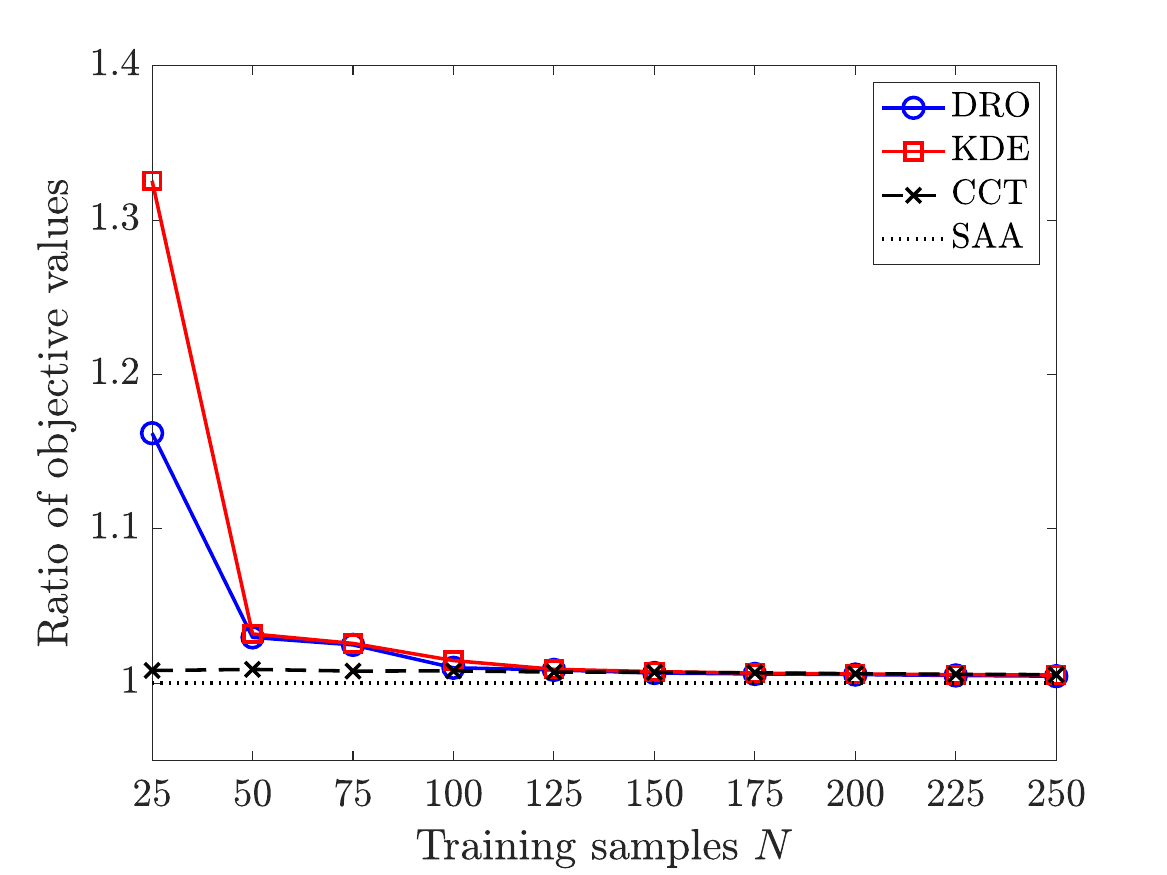}
\end{center}
\end{subfigure}%
\vspace{0.2cm}
\caption{{\textnormal{Probability of meeting the safety conditions (left) and transportation costs (right) for several data-driven approaches in our transportation problem with uniformly distributed demands. Both figures present median quantities over $100$ random instances.}} \label{fig:uniform}}
\end{figure}

\begin{figure}[tb]
\begin{subfigure}{.5\textwidth}
\begin{center}
\includegraphics[width=1\linewidth]{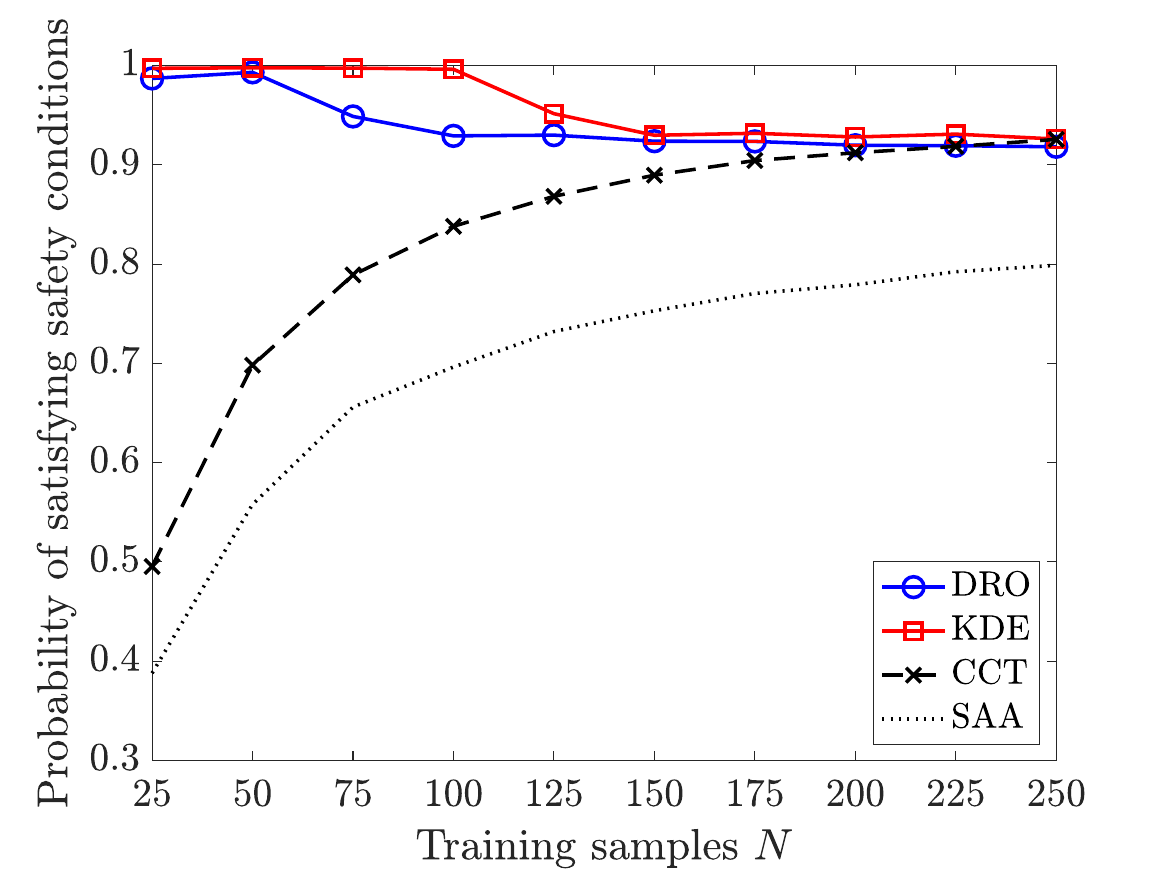}
\end{center}
\end{subfigure}%
\begin{subfigure}{.5\textwidth}
\begin{center}
\includegraphics[width=1\linewidth]{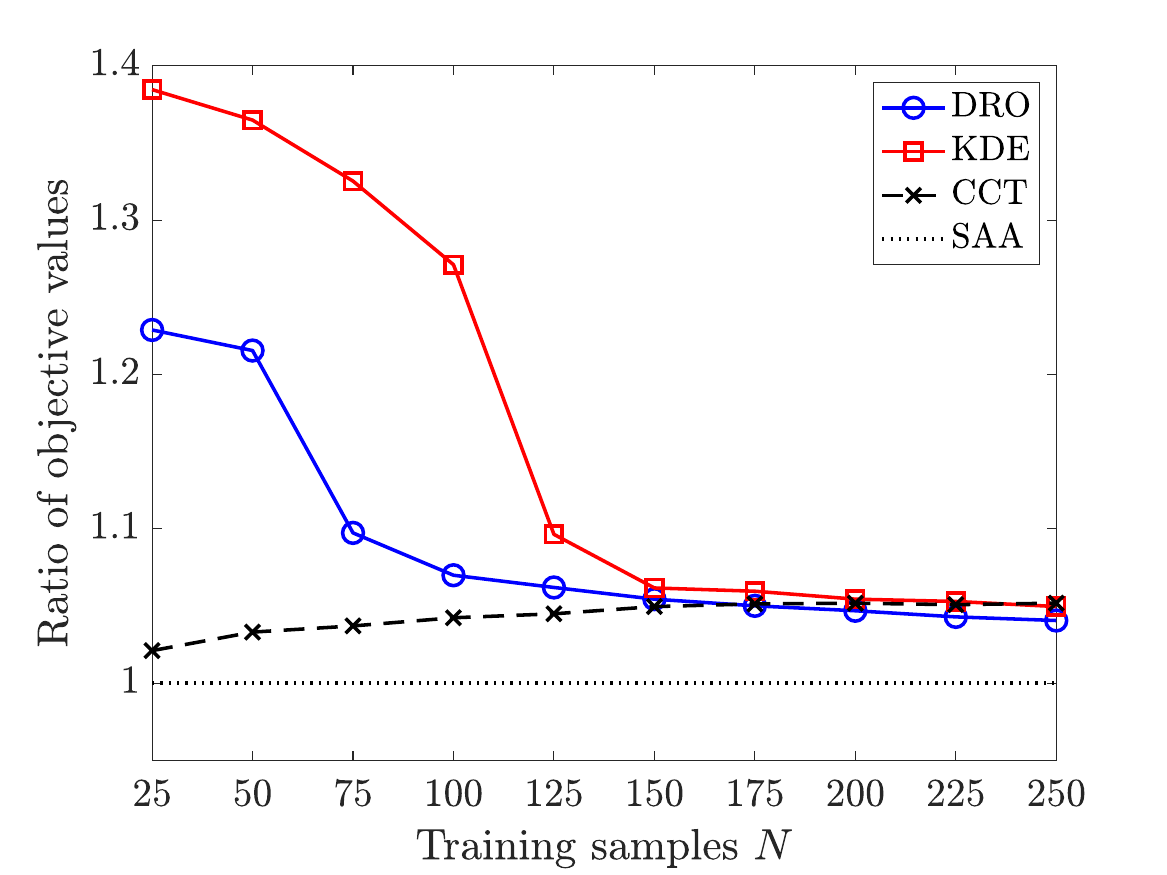}
\end{center}
\end{subfigure}%
\vspace{0.2cm}
\caption{{\textnormal{Probability of meeting the safety conditions (left) and transportation costs (right) for several data-driven approaches in our transportation problem with normally distributed demands. Both figures present median quantities over $100$ random instances.}} \label{fig:normal_unbounded}}
\end{figure}

\begin{figure}[tb]
\begin{subfigure}{.5\textwidth}
\begin{center}
\includegraphics[width=1\linewidth]{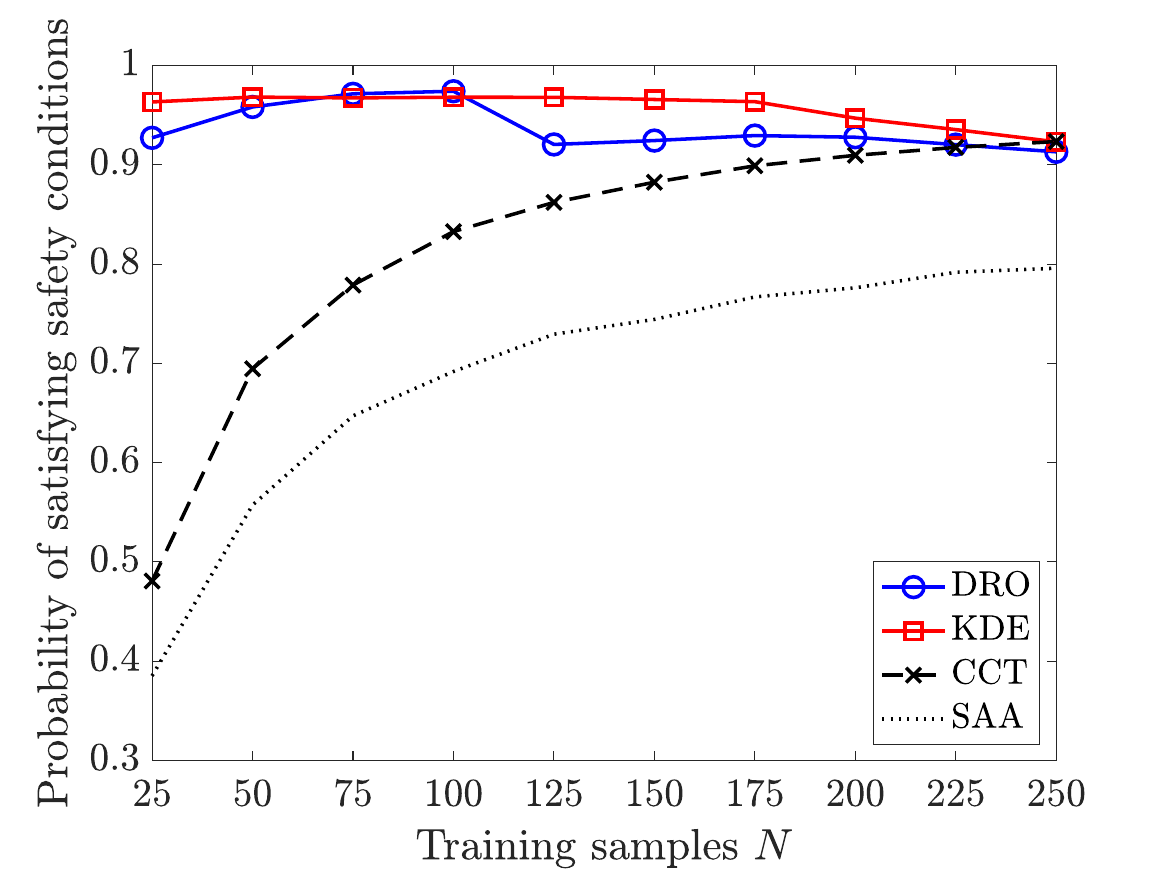}
\end{center}
\end{subfigure}%
\begin{subfigure}{.5\textwidth}
\begin{center}
\includegraphics[width=1\linewidth]{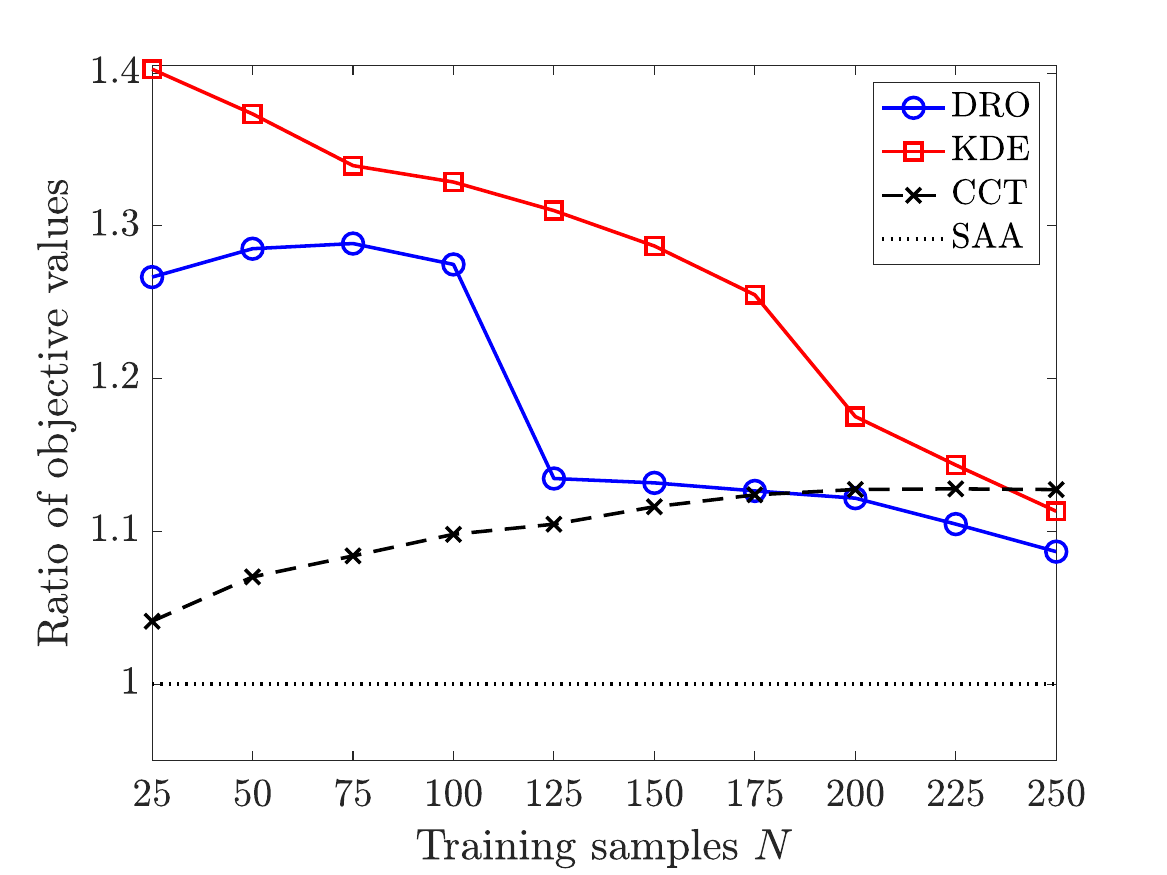}
\end{center}
\end{subfigure}%
\vspace{0.2cm}
\caption{{\textnormal{Probability of meeting the safety conditions (left) and transportation costs (right) for several data-driven approaches in our transportation problem with exponentially distributed demands. Both figures present median quantities over $100$ random instances.}} \label{fig:exponential_unbounded}}
\end{figure}

We next compare the out-of-sample performance of  our ambiguous chance constrained program~\eqref{eq:amb_transp_prob}, where the risk threshold $\varepsilon \in \{ 0.1, \, 0.05, \, 0.01 \}$ and the Wasserstein radius $\theta \in \{ 1\mathrm{E}-i \, : \, i = 2, 3, \ldots, 6 \}$ are selected using a $7$-fold cross-validation on the training dataset (`DRO'), with \emph{(i)} the classical chance constrained program~\eqref{eq:classical_transp_prob}, where the risk threshold is fixed to $\varepsilon = 0.1$ (`SAA'), \emph{(ii)} a variant of the classical chance constrained program~\eqref{eq:classical_transp_prob}, where the risk threshold $\varepsilon \in \{ 1\mathrm{E}-i \, : \, i = 1, 2, \ldots, 5 \} \cup \{ 0.05 \}$ is selected using a $7$-fold cross-validation on the training dataset (`CCT'), as well as \emph{(iii)} a Kernel density estimator based version of the ambiguous chance constrained program over a $\phi$-divergence ambiguity set, where the risk threshold $\varepsilon \in \{ 0.1, \, 0.05, \, 0.01 \}$ and the bandwidth $h \in \{ 1\mathrm{E}-i \, : \, i = -2, -1, \ldots, 3 \}$ of the Gaussian kernel are selected using a $7$-fold cross-validation on the training dataset (`KDE'; see \citealp{Jiang_Guan_2016}). We note that CCT can be regarded as a cross-validated version of the `best data-driven reformulation' proposed by \citet{lam:best}. We generate random problem instances with $5$ factories, $20$ distribution centers and $25$, $30$, \ldots, $250$ training samples. In all experiments, the expected demand $\mu_d$ at distribution center $d \in [D]$ is picked uniformly at random from the interval $[0, 10]$, whereas the actual demands follow a uniform distribution that is supported on $[0.8 \bm{\mu}, 1.2 \bm{\mu}]$ (Figure~\ref{fig:uniform}), a normal distribution with mean $\bm{\mu}$ and covariance matrix $0.1 \cdot \text{diag} (\bm{\mu})$ (Figure~\ref{fig:normal_unbounded}) or an exponential distribution where each distribution center $d \in [D]$ faces a demand $(1 + 0.4 \cdot [\tilde{\zeta}_d - 0.5]) \cdot \mu_d$, where $\tilde{\zeta}_d$ follows an exponential distribution with parameter $\lambda = 2$ (Figure~\ref{fig:exponential_unbounded}). In all cases, the demands are truncated to the non-negative real line. Our results indicate that the classical chance constrained program~\eqref{eq:classical_transp_prob} generates solutions that significantly violate the chance constraint, even if we select the risk threshold $\varepsilon$ out-of-sample. The two ambiguous chance constrained formulations, on the other hand, achieve the desired risk threshold, often at a modest increase in transportation costs. While our approach and the $\phi$-divergence ambiguity set perform similarly, our formulation appears to result in lower transportation costs, especially when data is scarce.

\section*{Acknowledgments}
The authors are grateful to the review team for constructive comments that led to substantial improvements of the paper.
The authors gratefully acknowledge financial support from the ECS grant~9048191, the SNSF grant BSCGI0$\underline{~}$157733 and the EPSRC grant EP/N020030/1.

\bibliographystyle{ormsv080}
\bibliography{References}

\newpage
\begin{appendices}
\section{Distance to a Union of Halfspaces}		
The distance of a point $\bmh{\xi} \in \mathbb{R}^K $ to a closed set $\mathcal{C} \subseteq \mathbb{R}^K$ with respect to a norm $\|\cdot\|$ is defined as 
$$
\mathbf{dist}(\bmh{\xi}, \mathcal{C}) = \min\{\|\bm{\xi} - \bmh{\xi}\| \mid \bm{\xi} \in \mathcal{C}\}.
$$
Note that the minimum is always attained. In the following, we derive a closed-form expression for the distance of a point to the union of finitely many closed halfspaces.
\begin{lemma}
\label{lem:distance to the union of closed half-spaces}
Let $\mathcal{H}_m = \{\bm{\xi} \in \mathbb{R}^K \mid a_m \geq \bm{b}^\top_m \bm{\xi}\}$ be a closed halfspace for each $m \in [M]$. If $\mathcal{C} = \bigcup_{m \in [M]} \mathcal{H}_m$ denotes the union of all halfspaces, then the distance of a point $\bmh{\xi}$ to $\mathcal{C}$ is given by
$$
\mathbf{dist}(\bmh{\xi}, \mathcal{C}) = \min_{m \in [M]} \bigg\{\dfrac{(\bm{b}^\top_m\bmh{\xi} - a_m)^+}{\|\bm{b}_m\|_*}\bigg\} = \bigg(\min_{m \in [M]} \bigg\{\dfrac{\bm{b}^\top_m\bmh{\xi} - a_m}{\|\bm{b}_m\|_*}\bigg\}\bigg)^+.
$$
\end{lemma}
\noindent \emph{Proof of Lemma~\ref{lem:distance to the union of closed half-spaces}.} $\;$
We first prove the assertion for $M=1$, in which case $\mathcal{C} = \mathcal{H}_1$. 
We thus have
\begin{align*}
\mathbf{dist}(\bmh{\xi}, \mathcal{C}) =~& \min_{\zeta, \bm{\xi}} \big\{ \zeta \mid \zeta \geq \|\bm{\xi} - \bmh{\xi}\| , ~ a_1 \geq \bm{b}^\top_1\bm{\xi} \big\} \\
=~ &  \max_{u, \bm{v}, w} \big\{\bm{v}^\top\bmh{\xi} - w a_1  ~\big|~  u = 1 ,~ \bm{v} =  \bm{b}_1 w, ~ u \geq \|\bm{v}\|_*, ~w \geq 0 \big\}\\
=~ &  \max_{ w} \big\{(\bm b_1^\top \bmh{\xi}- a_1) w ~\big|~ w\le 1/ \|\bm{b}_1\|_*, ~w \geq 0\big\} \\
=~ & \dfrac{(\bm{b}^\top_1\bmh{\xi} - a_1)^+}{\|\bm{b}_1\|_*},
\end{align*}
where the second equality follows from strong conic duality, which holds because the primal minimization problem is strictly feasible. Similarly, for $M\ge 1$ we find
\begin{align*}
\mathbf{dist}(\bmh{\xi}, \mathcal{C}) =~& \min_{m \in [M]} \mathbf{dist}(\bmh{\xi}, \mathcal{H}_m)= \min_{m \in [M]} \bigg\{\dfrac{(\bm{b}^\top_m\bmh{\xi} - a_m)^+}{\|\bm{b}_m\|_*}\bigg\} = \bigg(\min_{m \in [M]} \bigg\{\dfrac{\bm{b}^\top_m\bmh{\xi} - a_m}{\|\bm{b}_m\|_*}\bigg\}\bigg)^+,
\end{align*}
where the second equality follows from the first part of the proof. 
\hfill \Halmos
\endproof

\end{appendices}
\end{document}